\documentclass[times,11pt]{article} 
\usepackage{amssymb} 
\usepackage{amsmath} 
\usepackage{latexsym} 
\usepackage{epsf} 
\usepackage[all]{xy}  
\usepackage{proof} 
\usepackage{lscape} 
\usepackage{diagrams}  
\usepackage{verbatim}    
\usepackage{times}  
 
\newcommand{\rem}[1]{\relax} 
\newcommand{\cut}[1]{\relax} 

\newcommand{\semantics}[1]{[\![ #1 ]\!]}

\renewcommand{\phi}{\varphi}

\newtheorem{Th}{Theorem}[section]  
\newtheorem{theorem}[Th]{Theorem}  
\newtheorem{proposition}[Th]{Proposition}

\newtheorem{definition}[Th]{Definition}

\newcommand{\bit}{\begin{itemize}} 
\newcommand{\eit}{\end{itemize}\par\noindent} 
\newcommand{\ben}{\begin{enumerate}} 
\newcommand{\een}{\end{enumerate}\par\noindent} 
\newcommand{\beq}{\begin{equation}} 
\newcommand{\eeq}{\end{equation}\par\noindent} 
\newcommand{\beqa}{\begin{eqnarray*}} 
\newcommand{\eeqa}{\end{eqnarray*}\par\noindent}  
\newcommand{\beqn}{\begin{eqnarray}}   
\newcommand{\eeqn}{\end{eqnarray}\par\noindent}

\def\endproof{\hfill$\Box$}  
\def\qed{\ifmmode 
          $\Box$ 
         \else 
         {\unskip 
          \nobreak 
          \hfil 
          \penalty50 
          \hskip1em 
          \null 
          \nobreak 
          \hfil 
          $\Box$ 
          \parfillskip=0pt 
          \finalhyphendemerits=0 
          \endgraf} 
         \fi} 

\def\suck{\vspace{-2mm}} 
 
\title{\bf Epistemic actions as resources}  
\author{} 
\date{}  
 
\begin{document}  
\maketitle 
  
\vspace{-1.4cm}  
\noindent 
\begin{center} 
\begin{tabular}{cc} 
Alexandru Baltag {\ \ }and{\ \ } Bob Coecke &  
Mehrnoosh Sadrzadeh\\ 
\hspace{0.5cm}Oxford University Computing laboratory\hspace{0.5cm} &  
University of Southampton\\ 
{\sf\footnotesize baltag / coecke@comlab.ox.ac.uk} & 
\hspace{0.5cm}{\sf\footnotesize ms6@ecs.soton.ac.uk}\hspace{0.5cm} 
\end{tabular} 
\end{center} 
  
\vspace{2mm}\noindent 
\begin{abstract}  
\noindent  
We provide \emph{algebraic semantics}  together with a sound and complete \emph{sequent calculus}  for \em information update \em due to \em epistemic actions\em. This semantics is flexible enough to accommodate  incomplete as well as wrong information e.g.~due to \em secrecy \em and \em deceit\em, as well as \em nested knowledge\em. We give a purely algebraic treatment of the muddy children puzzle, which moreover extends to situations where the children are allowed to lie and cheat.  Epistemic actions, that is, information exchanges between agents $A,B, \ldots \in {\cal A}$, are modeled as elements of a \emph{quantale}.
The quantale $(Q,\bigvee,\bullet)$ acts on an underlying \emph{$Q$-right module} 
$(M,\bigvee)$ of  epistemic propositions and facts.  The epistemic 
content is encoded by \emph{appearance maps}, one pair 
$f^M_A\colon M\to M$ and $f^Q_A\colon Q\to Q$ of (lax) morphisms for each agent 
$A\in{\cal A}$, which preserve the module and quantale structure respectively. 
By adjunction, they give rise to \emph{epistemic modalities}~\cite{FH}, capturing the agents' \emph{knowledge}  on propositions and actions.  The module action 
is \emph{epistemic update}  and gives rise to 
\emph{dynamic modalities}~\cite{Kozen} --- cf.~\em weakest 
precondition\em. This model subsumes the crucial fragment of Baltag, Moss and Solecki's  \cite{BaltagMossSolecki} \emph{dynamic epistemic logic}, abstracting it in a constructive fashion while introducing \em resource-sensitive \em structure on the epistemic actions.  
\end{abstract} 
 
\section{Introduction} 
 
Consider the following well-known  puzzle.  After $n$ children played in the mud $k$ of them 
have mud on their forehead. They can see each other's foreheads but not their own ones. Their 
father initially announces ``At least one of you has mud on his forehead!''.  Then he asks: ``Is it you who has mud on his forehead?''.   Typically the children will all together answer: ``I don't know!''.
Again father asks: ``Is it you who has mud on his forehead?'', and again typically the children will all together answer: ``I don't know!''.  It turns out that after $k-1$ rounds of father's question and 
the children's ``I don't know!''-answers the ones which have mud on their forehead 
will now all know this.  Indeed, in the case of $k=1$ the dirty child knows that it must be him who is dirty since all the other children are clean.  In the case $k=2$ the two dirty children see only one other dirty child, so after a round of ``I don't know!''-answers, they realise that they must be dirty since in the case of only one dirty child that child should have known this already in the first round.  This argument extends to arbitrary $k>2$ by induction. 

This \em Muddy Children Puzzle \em exposes the need for a logical account of actions and agents as \em dynamic \em and \em epistemic resources \em  in situations involving 
information  exchange. Indeed, repetition of the same announcement provides new information to 
the children.    In particular, these \em dynamic  resources \em constitute the \em use-only-once resources \em of Girard's \em Linear Logic \em \cite{Girard}. In  linear logic, as compared to ordinary logic,  premisses cannot be copied nor deleted.\footnote{Strictly speaking we are only considering the multiplicative fragment of linear logic of which the non-linear counterpart is intuitionistic logic.}  We will also deal with  \em epistemic resources\em: presence of agents within a context, or availability of these agents as computing resources for other agents, 
affects the validity of deductions and execution of some actions by other agents. 
E.g.~for the children to make the correct conclusion they need to take into account the capabilities of  other children to make deductions.  In other words, some deductions are only valid in the presence of certain other agents. Our intuitionistic  sequences  
\[  
m_1,\ldots,q_1,\ldots,A_1,\ldots,m_k,\ldots,q_l,\ldots,A_n\vdash\delta 
\] 
consist of  different types of formulas and for example can contain propositions $m_1,\ldots,m_k$, actions $q_1,\ldots,q_l$ and agents $A_1,\ldots A_n$,  which resolve into a single proposition or action $\delta$. Also, a deduction might not be valid in the ``real world" while it is valid in the world as it \em appears \em to an agent.

To cast all this in mathematical terms we rely on order-theory.  A proposition $a$ is implied by a proposition $b$ iff $b\leq a$. An action $p$ is less deterministic than an action $q$, i.e.~$p=q$ `or' $p=q'$, iff $q\leq p$.   Actions can be performed one after the other and in order to be able to reason with them we require distributivity of `composition' over `or'.  The resulting structure is a \em quantale\em.  Quantales  have been used as semantics for non-commutative Intuitionistic Linear Logic~\cite{Yetter}, which itself traces back to Lambek calculus~\cite{Lambek}.\footnote{Quantales are to complete Heyting algebras what monoidal closed categories are to Cartesian closed categories, respectively providing semantics for Intuitionistic Logic,  and for non-commutative Intuitionistic  Linear Logic, including Lambek calculus.}  The quantale acts on a  \em  sup-lattice \em of propositions.  Both the quantale and the sup-lattice come equipped with modal operators which capture the epistemics. 
They will allow to encode \em incomplete knowledge\em, e.g.~due to secrecy, \em wrong knowledge\em, e.g.~due to deceit, and \em nested knowledge\em, i.e.~one agent's knowledge on some other agent's knowledge, possibly yet again about some other agent's knowledge, and so on.
Technically these modal operators are so-called (lax-)endomorphisms of the above structure, one endomorphism-pair for each agent. Their Galois adjoints will stand for knowledge.  The pair of a sup-lattices and a quantale  without the modal operators have previously been used in concurrency  \cite{AbramskyVickers,Resende} and quantum logic \cite{CoeckeMS}.  Boolean algebras with adjoint operators, called Galois algebras, have previously been used in temporal logic~\cite{vonKarger}.   

Our  algebraic semantics and sound and complete sequent calculus further \em conceptualize \em and \em abstract \em the usual Kripke semantics  and Hilbert-style axiomatic logic for such situations e.g.~the dynamic epistemic logic of Baltag, Moss and Solecki~\cite{BaltagMoss,BaltagMossSolecki} {\bf (BMS)},
which is a {\bf PDL}-style logic to reason about epistemic actions and updates in a multi-agent system. 
Applications are secure communication, where issues of privacy, secrecy and authentication of  
communication protocols are central, software reliability for concurrent programs, AI,
where agents are to be provided with reliable tools to reason 
about their environment and each other's knowledge,
e-commerce, where agents need to have knowledge  
acquisition strategies over complex networks. The standard approach to information flow in a 
multi-agent system has been presented in~\cite{FH} but it does not present a 
formal description of epistemic actions and their updates. The first attempts to 
formalize such actions and updates were done by Gerbrandy and Groenveld \cite{Ger1,Ger2, GerGroe} and Plaza \cite{Plaza}, but they only studied a restricted class of these actions. A general notion of epistemic actions and updates was introduced in 
\cite{BaltagMoss,BaltagMossSolecki}. However, in this approach there is no account of resources in the underlying logic, and more importantly, the operations of sequential composition of actions and updating are concrete constructions on Kripke structures, rather than being
taken as the fundamental operations of an abstract algebraic signature.   In view of the purely Boolean nature of these Kripke models it is also worth stressing that in our proof of the Muddy children puzzle we essentially only \em reason by adjunction\em, both in terms of dynamic and epistemic residuals, but \underline{not} assuming the lattice of proposition to have complements nor for it to be distributive.
 
 We proceed as follows. First we introduce the objects of our algebra,  \em epistemic systems\em,   and justify their axiomatic 
structure.  We use our setting to analyze   
the Muddy Children Puzzle and  some of its more interesting and newer variants, involving lying children or secret communication,
as well as (a simplified version of) the Man-In-The-Middle (MITM) cryptographic attack.  We give examples of our structure and briefly explain how   models of  {\bf BMS}~\cite{BaltagMossSolecki} are instances of it, referring the reader for details of construction to~\cite{SadrThesis}. Next we introduce 
the sequent calculus and give a summary  of the completeness proof, referring the reader for the full  proof to~\cite{SadrThesis}.  We illustrate the use of the sequent calculus by proving a weak permutation property for  our epistemic and dynamic modalities and by encoding and deriving a property of  the MITM attack. We   conclude with suggestions for further elaboration. 
 
\section{The algebra of epistemic actions and epistemic propositions}  
 
A {\it sup-lattice} $L$ is a  complete lattice and a {\it sup-homomorphism} is a map between sup-lattices which preserves arbitrary joins.  We denote the \em bottom \em 
and \em top \em of $L$ by $\bot$ and $\top$ respectively, and its \em atoms \em by $Atm(L)$. 
A sup-lattice  is \em atomistic \em iff each element can be written as the supremum of the atoms below 
it.  Every sup-homomorphism $f^*:L\to M$ has a right Galois adjoint $f_*:M\to L$, i.e.
\[
f^*(a)\leq b\Leftrightarrow a\leq f_*(b)\,,
\]
which preserves arbitrary infima. We denote an adjoint pair by  
$f^*\dashv f_*$.  In computational terms, the right Galois adjoint $f_*$ assigns \em weakest 
preconditions \em to its arguments, given the \emph{program} $f^*$.   
 
A {\it quantale} is a sup-lattice $Q$ with a monoid structure 
$(Q, \bullet, 1)$  which distributes over arbitrary joins at both sides.
Since for all $a\in Q$ the maps $a\bullet -:Q\to Q$ and $-\bullet 
a:Q\to Q$ preserve arbitrary joins they have right Galois adjoints 

\[
a \bullet - \,\dashv\, a \setminus - \qquad\qquad{\rm and}\qquad\qquad - \bullet\ a \,\dashv\,  -\, /\, a\,,
\]
explicitly given by 
\[ 
a \setminus b := \bigvee \{c \in Q \mid a \bullet c 
\leq b\}\qquad{\rm and}\qquad 
b\,/\,a := \bigvee\{c \in Q \mid c \bullet a \leq b\}. 
\] 
A map $f:Q\to Q$ is a \em quantale homomorphism \em if it is both a sup-homomorphism and a 
monoid-homomorphism. It is a \em lax quantale homomorphism \em if it is a sup homomorphism 
and if 
\[
1\leq f(1)\qquad\qquad{\rm and}\qquad\qquad f(a\bullet b)\leq f(a)\bullet f(b)\,.
\]

Examples of quantales are: the set ${\sf sup}(L)$ of all 
sup-endomorphisms of a complete lattice $L$ ordered pointwisely; the set of 
all relations from a set $X$ to itself ordered by pointwise inclusion --- this 
quantale is isomorphic to ${\sf sup}({\cal P}(X))$; the powerset of any 
monoid with composition extended by continuity. 

Since quantales are monoidal closed categories they provide a semantics for non-commutative Intuitionistic Linear Logic~\cite{Yetter,Girard,AbramskyVickers}: linearity of monoidal closed categories 
follows by the 
\em absence 
\em (in general) of natural morphisms $\Delta_A:A\to A \otimes A$ and left and right projections $p_1:A\otimes B\to 
A$ and $p_2:B\otimes A\to A$, and hence quantales (in general) do not satisfy $\,a\leq a\bullet a\,$ nor $\,a\bullet b\leq a\,$ nor $\,a\bullet b\leq b\,$ (where now $\bullet$ is the monoidal tensor $\otimes$).   Note that quantales have more operators (than multiplicatives), with regard to which they are not resource-sensitive, for example we have similar inequalities for the meet of the quantale, that is we have that $a \leq a \wedge a$, and also $\,a \wedge b\leq a\,$ and $\,a\wedge b \leq b$.  
   
A \em $Q$-right module \em for a quantale $Q$ is a 
sup-lattice $M$ with a \em module action \em 
\[
-\cdot-: M \times Q \rightarrow M
\]
which preserves arbitrary joins in both arguments, 
\[
m \cdot 1= m\qquad\qquad {\rm and}\qquad\qquad
m \cdot (q_1 \bullet q_2) = (m \cdot q_1) \cdot q_2
\]
for all $m\in M$ and all $q_1,q_2\in Q$.  We have two adjoint pairs
$- \cdot q \dashv [q]-$ and $m \cdot - \dashv \{m\}-$ 
where 
\[ 
 [q]m := \bigvee\{m' \in M \mid m' \cdot q \leq m\}\qquad{\rm and}\qquad
\{m\}m' := \bigvee \{q \in Q \mid m \cdot q \leq m'\}. 
\] 

For example, a quantale $Q$ is a $Q$-right module over itself with 
composition as the action and a complete lattice $L$ is a ${\sf 
sup}(L)$-right module with function application as the action.  
For details on quantales, $Q$-modules and also 
$Q$-enrichment we refer the reader to \cite{JoyalT,Paseka,Rosenthal,Stubbe}.  
For applications of these in computing, linguistics and physics we refer the reader to 
\cite{AbramskyVickers,CoeckeMS,Hoare,Lambek,Mulvey,Resende}.  
  
\begin{definition}\em  \cite{AbramskyVickers}
A {\it system} is a pair $(M,Q)$  
with 
$Q$  a quantale and $M$ a $Q$-right module.  
\end{definition}  

\begin{definition}\em 
A \em system-endomorphism \em  
$(M,Q) \rTo^{f} (M,Q)$ 
is a pair $\left(f^{M}: M \rightarrow M\,,\, f^{Q}: Q \rightarrow Q\right)$ 
where $f^{M}$ and $f^{Q}$ are both  sup-homomorphisms,  and for all $m\in M$ and $q, q'\in 
Q$ we have 
\begin{equation}\label{eq:mult}
f^Q (q\bullet q')\, \leq\, f^Q (q) \bullet f^Q (q')
\end{equation} 
\begin{equation}\label{eq:update} 
\,f^{M}(m \cdot q)\, \leq\, f^{M}(m) \cdot f^{Q}(q)
\end{equation} 
\vspace{-4mm}
\begin{equation}\label{eq:skip}
\,\ 1\, \leq\, f^Q (1)\,.
\end{equation} 
\end{definition}  
Hence $f^Q$  is \em lax functorial \em and a lax quantale homomorphism\footnote{Our notion of system endomorphism also differs from the one in the literature, (e.g.~\cite{JoyalT} and categories of modules for rings) in that we consider non-trivial homomorphisms on the quantale, a so-called \em change of base\em.  Explicitly, we do not have $f(m \cdot q) = f(m) \cdot q$ \ for a system endomorphism $f$.}. 

\begin{definition}\label{def:epistemicsystem}\em 
An \em epistemic system \em is a tuple $(M, Q, \{f_{A}\}_{A \in {\cal A}})$ where $(M,Q)$ 
is a  system and 
$\{f_{A}\}_{A \in {\cal 
A}}$ are system-endomorphisms. The elements of ${\cal A}$ are called \em agents\em, the elements of 
$Q$  \em epistemic actions \em and the elements of $M$  \em epistemic propositions\em. The system-endomorphisms are called \em appearance maps\em.
\end{definition}  
 
 
\paragraph{Epistemic Propositions.} 
We interpret the elements of the module as \emph{epistemic propositions} and  their order relation $m\leq m'$ for  
$m,m'\in M$ as \em logical entailment \em $m\vdash m'$. The epistemic proposition 
$f_A^M(m)$ describes how the world appears to agent $A$: it comprises all propositions that agent $A$ 
\em believes to hold \em whenever $m$ holds in the `real world'. Two extreme examples are $f_A^M(m)=\top$, which corresponds to absence of any 
knowledge whatsoever, and 
$f_A^M(m)=m$, which stands for complete knowledge. If for $m,m'\in M$ we have $f_A^M(m)<f_A^M(m')$ then agent $A$ possesses strictly more (possibly incorrect) knowledge on $m$ than on  $m'$.  It also follows that $f_A^M$ indeed needs to be covariantly monotone 
--- the additional preservation of suprema will assure existence of epistemic modalities (see below).  If for $m\in M$ we have $f_A^M(m)<f_B^M(m)$ agent $A$ possesses strictly 
more (possibly incorrect) knowledge on $m$ than agent $B$.  But indeed, this knowledge is not necessarily correct! If 
for $m, m' \in M$ with $m\not\leq m'$ we have 
$f_A^M(m)\leq m'$ then agent $A$ \emph{believes incorrect information} to be true, e.g.~due to 
deceit of another agent, a malfunctioning communication channel, corrupted data etc.    
If the module is atomistic, then the atoms can be thought of  as states --- cf.~Kripke structures representing epistemic scenarios (see also the following section and \cite{BCS} and  \cite{SadrThesis}). 
 
\paragraph{Knowledge {\rm and/or} belief.} 
For each agent $A\in{\cal A}$ let $\square_{A}^{M}$ be defined by~$f_{A}^{M}\dashv \square_{A}^{M}$. By adjunction we have $m\leq\square_{A}^{M}m'$ if and only if $f_{A}^{M}(m)\leq m'$, that is, ``when proposition $m$ holds, agent  $A$ knows/believes $m'$''.  Hence $\square_{A}^{M}m$ stands for \em agent $A$'s knowledge/belief \em on $m$. 
This modality indeed covers both knowledge and belief: in contexts where no wrong belief is allowed, we read it as knowledge or \em justified true  belief\em, and otherwise, as \em justified belief\em. Since $\square_{A}^{M}$ is a right Galois adjoint we have $\square_{A}^{M}(\bigwedge_im_i)=\bigwedge_i\square_{A}^{M}m_i$.
Hence it preserves the empty and binary meets, and is monotone: 
\[ 
\square_{A}^{M}\top=\top\qquad\qquad\qquad\square_{A}^{M}(m\wedge   
m')=\square_{A}^{M}m\wedge\square_{A}^{M}m' 
\qquad\qquad\qquad 
{m\leq m' 
\over 
\square_{A}^{M}m\leq \square_{A}^{M}m'}\,.  
\]  
When $M$ is a \em frame \em (= complete Heyting algebra  \cite{Johnstone}) we can internalize the partial order using the defining 
property of a Heyting algebra. 
In the special case that $Q=\{1\}$ and $A=\{*\}$ we obtain    
the intuitionistic modal logic ${\bf IntK}_{\square}$ of 
\cite{WolterZakhatyaschev}.  
If $M$ is moreover a complete 
boolean algebra (e.g.~the powerset of its atoms) then Kripke's axiom {\bf K} follows i.e. 
\[ 
\square_{A}^{M}(m\to m')\to 
(\square_{A}^{M}m\to \square_{A}^{M}m').  
\] 
Diamonds and corresponding rules arise in that case by duality.   
If $M$ is atomistic and the set of atoms are denoted by $S$ then to each $f_A^M$ one can assign an \em accessibility relation{\,} \em $\stackrel{A}{\,\to\ \,}\subseteq S\times S$ by setting  \[
s\stackrel{A}{\to}s' \Longleftrightarrow s'\leq  f_A^M(s)\,. 
\]
It is this relation which is primitive in  ordinary epistemic logics rather than appearance maps. 
But in our setting, in general, this accessibility relation turns out not to be reflexive, nor (anti-)symmetric, nor transitive e.g.~\em positive introspection \em  ${\square_{A}^{M}m\leq\square_{A}^{M}\square_{A}^{M}m}$ does not hold in general. 
 
\paragraph{Epistemic actions.} 
We interprete the elements of the quantale as \emph{epistemic actions}  where the order is \em information ordering\em: if for  $q,q'\in Q$ we have  $q\leq q'$ then $q'$ is less deterministic than $q$. 
The suprema $\bigvee_iq_i$ in the  
quantale, similar as in~\cite{AbramskyVickers,CoeckeMS}, stand for \em non-deterministic choice\em. The action $f_A^Q(q)$ captures how  $q$ appears to agent $A$. The appearance maps allow  to accommodate actions such as \em information hiding \em or \em encryption, \em by  ${q<f_A^Q(q)}$, and 
\em misinformation \em such as lying, cheating and deceit by $q\not\leq f_A^Q(q)$. 
Analogously to the case of propositions, setting $f_{A}^{Q}\dashv \square_{A}^{Q}$
stands for \em agent $A$'s knowledge/belief on \em $q$ i.e.~``when action $q$ is happening, agent $A$ believes action $q'$ to be happening''. 
These epistemic modalities $\square_{A}^{Q}$ satisfy the same properties as $\square_{A}^{M}$. 
If the quantale is atomistic then its atoms can be interpreted as \em deterministic \em actions. 

\paragraph{Sequential composition.} 
The quantale multiplication stands for sequential composition of epistemic actions.  
The multiplicative unit $1$ is the \em void \em epistemic action, that is, \em nothing happens\em, sometimes referred to as \em skip \em in literature (cf.~~\cite{Kozen}).  We do not require $f_A^Q (1)=1$ but only $1\leq f_A^Q (1)$ since this enables us to accommodate \em suspicions\em, cf.~eq.(\ref{eq:skip}).  By this we mean that even when nothing is happening one could still suspect that something hidden might be happening, say $q$, resulting in $f_A^Q (1)=1\vee q$.  Suspicions are for example important for applications to protocol security, see~\cite{SadrThesis} ch. 5 for some examples. On the other hand requiring $1\leq f_A^Q (1)$ imposes \em rationality \em of the agent (vs.~insanity): if nothing is happening then the agent considers nothing to be happening at least as an option.  This argument carries over to \em appearance of sequential composition\em, again subject to a rationality requirement, and suspicions cause \em laxity\em, cf.~eq.(\ref{eq:mult}):
\[
f^Q_A(q\bullet 1)=f^Q_A(q)=f^Q_A(q)\bullet 1\leq f_A^Q (q)\bullet f^Q_A(1)\,.
\]
Other situations where we have a strict inequality arise when $q \bullet q' = \bot$ and thus $f_A^Q (q \bullet q') = \bot$, but $f_A^Q (q) \bullet f_A^Q (q') \neq \bot$ again due to the fact that the agent might suspect more options then what is actually happening --- for a detailed discussion and a concrete example see~\cite{SadrThesis}. 
 
\paragraph{Epistemic updating.} 
The action of the quantale on the module encodes the crucial notion of \emph{epistemic updating}.
After performing an epistemic action 
$q\in Q$ on an epistemic proposition $m\in M$ we obtain a new epistemic proposition $m\cdot q\in M$.  
Each agent updates his knowledge according to how he perceives the epistemic action, so $f_{A}^M(m\cdot  q)$ relates to $f_{A}^M(m) \cdot  f_{A}^Q(q)$. Again suspicions impose laxity,
cf.~eq.(\ref{eq:update}):
\[
f^M_A(m\cdot 1)=f^M_A(m)=f^M_A(m) \cdot 1\leq f^M_A(m)\cdot f_A^Q (1)\,,
\]
and we can have situations where  an action $q$ cannot apply to a proposition $m$, that is $ m \cdot  q = \bot$, and thus $f_A^M (m \cdot  q) = \bot$, but the appearance of the action can apply to the appearance of the proposition, that is $f_A^M(m) \cdot  f_A^Q(q) \neq \bot$ --- for a detailed discussion see~\cite{SadrThesis}. Situations where some of the suspected alternatives yield contradiction after update yield a process of \emph{learning} (or acquiring more information): the agent will eliminate his contradiction-leading views  and not anymore consider them as true options. 
 
\paragraph{Dynamic modalities.}  
Since both update $-\cdot-$ and quantale multiplication $-\bullet-$ preserve suprema in both arguments, a 
range of {\it residuals} arise, namely  
\[ 
- \cdot  q \dashv [q]-  
\quad\quad\quad\quad  
m \cdot  - \dashv \{m\}- 
\quad\quad\quad\quad  
q \bullet - \dashv q \setminus - 
\quad\quad\quad\quad  
- \bullet \, q \dashv  - / q 
\] 
for each $m\in M$ and each $q\in Q$. The residual $[q]-$ is the dynamic modality of  dynamic logic 
\cite{Kozen}, that is, \em weakest precondition\em. We read $[q]m$ as ``after program $q$ proposition $m$ holds''. On the other hand, $m\cdot q$ is the \emph{strongest postcondition}. The other ones are variants on these e.g. see~\cite{Hoare}. In particular the ones with respect 
to sequential composition correspond to the \em residuals 
\em of Lambek calculus \cite{Lambek} and the linear implications of non-commutative Linear Logic.   
  
\paragraph{Kernel.} 
If $m\cdot  q=\bot$ then $q$ cannot be applied to $m$. We define a  \em kernel 
\em for an action 
$q\in Q$ as
\[ 
Ker(q) := \{m \in M \mid m \cdot  q = \bot\},  
\] 
i.e.~as the \emph{co-precondition} of an action $q$ (= the dual to the so-called \em precondition \em of $q$). Since 
\[
Ker(q)=\downarrow\!\left(\bigvee Ker(q)\right)\,,
\]
 ``\em not \em being  in the precondition of $q$'' exists as a proposition in $M$ for 
all $q\in Q$. Also note that the kernel of each action is the weakest proposition to which  the action cannot  apply, that is $Ker(q) = [q]\bot$. 
 
\paragraph{Stable facts.} 
Each epistemic system has a non-epistemic part, referred to as  \emph{facts}, being the  propositions which cannot be 
altered by any epistemic action.  Define the \em stabilizer \em of $Q$ as
\[ 
Stab(Q) := \{\phi \in M \mid \forall q \in Q\,, \phi \cdot  q \leq \phi\}. 
\] 
It consists of those epistemic 
propositions which are stable under the epistemic actions, or equivalently, $\phi\leq[q]\phi$, which expresses preservation of validity of $\phi$: if it is true before running $q$, it will remain true afterwards. To summarize, epistemic 
propositions both encode actual facts and the knowledge of each agent, that is, they have both factual and epistemic 
content. 
 
\section{Examples of epistemic actions and epistemic systems}   
 
We present some examples of epistemic actions that can exist 
in an epistemic system $(M,Q,\{f_{A}\}_{A \in {\cal A}})$. 
\bit
\item {\bf Public refutation} 
of the proposition $m \in M$ is an epistemic action $q \in Q$ with $f_{A}^Q(q) = q$ 
for all 
$A\in{\cal A}$ and for which $Ker(q) = \downarrow\!m$. 
\item {\bf Private refutation to a subgroup} 
is an action that privately refutes 
$m$ to the subgroup $\beta$ of agents. In this case $Ker(q)$ is the 
same as above and $f_{A}^Q(q) = q$ for $A \in \beta$ and $f_{A}^Q(q)= 1$ otherwise. 
\item {\bf Failure test} 
of a proposition $m$ is an action $q$ 
that tests when $m$ fails. It is a particular case of private refutation 
where $m$ is refuted to an empty set of agents. Hence we have  $Ker(q)=\downarrow\!m$ and $f_A^Q(q) = 1$ for all 
$A\in{\cal A}$. 
\item {\bf Public announcement} is also definable 
in our setting. However, while ``being not in the precondition of $q$'' is a 
proposition in $M$ for all $q\in Q$, ``being in the precondition of $q$'' in general isn't one.  
To see this consider the lattice $\{\bot\leq a,b,c\leq 
\top\}$ with 
$q$ such that 
$Ker(q)=\{\bot,a\}$,  then both $b$ and $c$ are in the precondition but $b\vee c=\top$ isn't. The reason for 
this is that this lattice is non-Boolean with $a$ not having a complement. Hence public announcement of the 
proposition 
$m\in M$ is an epistemic action $q\in Q$ for which 
$f_{A}(q) = q$ and for which $\bigvee Ker(q)$ has a \em Boolean complement \em $(\bigvee Ker(q))^c$, 
satisfying 
$(\bigvee Ker(q))^c=m$. 
\item {\bf Private announcement to a subgroup} can be defined analogously. 
\eit

\paragraph{The Muddy Children Puzzle.}  
This  puzzle, explained in the introduction,  is a paradigmatic example in the standard epistemic logic literature --- e.g.~\cite{FH}.    In the usual encodings the communication between the father and children (i.e.~father's announcement and questions and the children's answers) is not part of the actual encoding. Our  approach (similar to the one in \cite{Baltag4}) does allow to encode communications and their effects on the agents'  knowledge. Our algebraic setting provides us, furthermore, with a semi-automatic elegant equational way of doing so. 

\smallskip
We encode the puzzle in an epistemic system. The set of agents ${\cal A}$ includes the children  
$C_1, \cdots, C_n$. We assume that $C_1, \cdots, C_k$ for $1\leq k\leq n$ are dirty. 
 The  module $M$ includes all possible initial propositions $s_\beta$ with $\beta\subseteq {\cal A}$ being those 
children that have mud on their forehead. 
For example $s_{C_1, \cdots, C_k}$ expresses the ``real state'' in which $C_1, \cdots, C_k$ are dirty and $C_{k+1}, \cdots, C_n$ are clean.  Since the children cannot see their own foreheads (which might  either be dirty or not) we have  
\[
f^M_{C_i}(s_\beta)=s_{\beta\setminus\{C_i\}}\vee s_{\beta\cup\{C_i\}}\,.
\] Let $D_\emptyset$ be the fact that no 
child has a dirty forehead and let $D_i$ be the fact that the $i$'th child has a dirty forehead, hence we have:
\[
\{D_\emptyset\}\cup\{D_i \in M \mid C_i \in {\cal A}\}\subseteq Stab(Q)\,.
\]
For the propositions and  facts we have $s_\beta\leq D_i$ for  all 
$C_i\in\beta$ and $s_\emptyset\leq D_\emptyset$, which sets that each proposition satisfies the corresponding fact.  Let $q\in Q$ be a round of  all children's  ``no'' answers i.e.~public refutation of $\bigvee_{i=1}^{i=n}\Box_{C_i}D_i$, hence 
$Ker(q)=\downarrow\!\bigvee_{i=1}^{i=n}\Box_{C_i}D_i$ and $f_{C_i}^Q(q) = q$ for $1\leq i\leq n$.   
Let $q_0\in Q$ be father's announcement that at least one child has mud on his forehead 
i.e.~$Ker(q_0)=\downarrow\!D_\emptyset$ and $f_{C_i}^Q(q_0) = q_0$ for $1\leq i\leq n$.  
 
\begin{proposition} 
After $k-1$ rounds of refutations, child $j$ for $1\leq j\leq k$ knows that he is dirty   
i.e. 
\begin{equation}\label{eq:muddychildren} 
s_{\{C_1, \cdots, C_k\}} \leq [q_0 \,(\bullet\, q)^{(k-1)}]\Box_{C_j}D_j
\end{equation} where $q_0 (\bullet\, q)^{(k-1)}$ denotes \ $q_0 \bullet q \bullet \cdots \bullet q$ with $k-1$ 
occurrences of $q$. 
\end{proposition} 
 
\noindent 
{\bf Proof.}  
We proceed by induction on the number $k$ of dirty children. 
If we move the dynamic modalities in eq.(\ref{eq:muddychildren}) to 
the left by adjunction we obtain  
\begin{equation}\label{eq:inducthyp} 
s_{\{C_1, \cdots, C_k\}} \cdot  q_0 \,(\cdot\  q)^{(k-1)}=s_{\{C_1, \cdots, C_k\}} 
\cdot  (q_0 \,(\bullet q)^{(k-1)}) \leq\Box_{C_j}D_j 
\end{equation}  
using the module structure.  
After moving the epistemic modality to the left and by the update inequality  eq.(\ref{eq:update}), it suffices to prove the following inequality
\[ 
f_{C_j}^M(s_{\{C_1, \cdots, C_k\}})\cdot  q_0 \,(\cdot\  q)^{(k-1)} \leq (s_{\{C_1, \cdots, C_k\}} \vee s_{\{C_1, \cdots, C_k\}\setminus\{C_j\}})\cdot  q_0 \,(\cdot\  q)^{(k-1)}\,
\]  
which is equivalent to the following by our assumption about  $f_{C_j}^M$
\[
(s_{\{C_1, \cdots, C_k\}} \vee s_{\{C_1, \cdots, C_k\}\setminus\{C_j\}})\cdot  q_0 \,(\cdot \ q)^{(k-1)}  \leq D_j\,. 
\]
By distributivity of $\vee$ over $\cdot$ and  the definition of suprema it suffices to 
prove  
\begin{equation}\label{eq:cases}  
s_{\{C_1, \cdots, C_k\}} \cdot  q_0  \,(\cdot \ q)^{(k-1)} \leq D_j 
\quad\quad {\rm and} \quad\quad 
s_{\{C_1, \cdots, C_k\}\setminus\{C_j\}} \cdot  q_0 \,(\cdot \ q)^{(k-1)} \leq D_j\,. 
\end{equation} 
We respectively refer to these inequalities as eq.(\ref{eq:cases}{\it l}\,) and eq.(\ref{eq:cases}{\it r}). 
First we show that eq.(\ref{eq:cases}{\it l}\,) holds for all $k$. 
Updating  both sides of $s_{\{c_1, \cdots, 
c_k\}}\leq D_j$ by $q_0 \,(\cdot \ q)^{(k-1)}$ we get 
\[s_{\{C_1, \cdots, C_k\}} \cdot  q_0 \,(\cdot \ q)^{(k-1)} \leq D_j \cdot  q_0 \,(\cdot \ q)^{(k-1)}\leq 
D_j 
\]  
where the last inequality follows by $D_j \in Stab(Q)$. Hence eq.(\ref{eq:cases}{\it l}\,).  
Now we prove the base case $k=1$ of our induction. 
Eq.(\ref{eq:cases}{\it r}) is $s_\emptyset 
\cdot  q_0\leq D_1$ in this case, which is true since $s_\emptyset\leq D_\emptyset\in Ker(q_0)$ so 
$s_\emptyset \cdot  q_0=\bot$. To prove 
eq.(\ref{eq:cases}{\it r}) we use the inductive hypothesis in terms of eq.(\ref{eq:inducthyp}).  
By symmetry of $\{C_1, \cdots, C_k\}$ we have\vspace{-2mm}  
\begin{equation}\label{eq:dynamics} 
 s_{\{C_1, \cdots, 
C_{k}\}\setminus\{C_j\}} \cdot  q_0 \,(\cdot \ q)^{(k-2)} \leq \Box_{C_j}D_j\leq 
\bigvee_{i=1}^{i=n}\Box_{C_i}D_i\vspace{-3.5mm} 
\end{equation} 
 
\noindent 
so 
\begin{equation}\label{eq:kernel}  
s_{\{C_1, \cdots, C_{k}\}\setminus\{C_j\}} \cdot  q_0 \,(\cdot \ q)^{(k-2)}\in Ker(q)
\end{equation}
 and hence 
\[ 
\bot=(s_{\{C_1, \cdots, C_{k}\}\setminus\{C_j\}} \cdot  q_0 \,(\cdot \ q)^{(k-2)})\cdot  q =s_{\{C_1, 
\cdots, C_{k}\}\setminus\{C_j\}} \cdot  q_0 \,(\cdot \ q)^{(k-1)}  \leq D_j 
\] 
i.e.~eq.(\ref{eq:cases}{\it r}), what completes the proof.  
\endproof\newline 

Analysing the dynamics of this proof we notice that in each inductive step we show that the epistemic state 
${s_{\{C_1, \cdots, C_k\}} \cdot  q_0 \,(\cdot \ q)^{(k-1)}}$ is included in the kernel of the refutation 
$q$ --- cf.~eq.(\ref{eq:dynamics}). This inductive update reflects the systematic update of the children's 
knowledge during the process.  Such a dynamics is not visible in the proofs performed in static epistemic 
logic \cite{FH} where there is no notion of update. 
 
This machinery not only enables us to deal with  classical epistemic scenarios in a dynamic way, but it also 
provides us with tools to treat (for the first time) other more complicated and realistic versions of these epistemic scenarios. 
As examples, we encode and analyze more complex versions of the above puzzle, in which some of the children may \em lie\em, 
or otherwise \em cheat \em by engaging in secret communication\footnote{The cheating example was done for Kripke models of {\bf BMS} by one of the authors~\cite{Baltag4}, while the lying example is new.}, as well as an example of a cryptographic attack. 

\paragraph{Lying Children.}
Assume that the same  $n$ children are playing in the mud and this time only one of them, say $C_1$, has a dirty forehead. Their father does the announcement exactly as in the classical  
Muddy Children Puzzle, and then asks the same question. Now before the first round of answers,  
the dirty child  who is a perfect reasoner, follows the proof presented above and by looking around and seeing no 
other dirty child, concludes that he is dirty $\Box_{C_1} D_1$. But instead of announcing the truth in the first round, 
he lies by saying that he does not know that he is dirty.  
This version is encoded using the same epistemic system as muddy children above with the difference that this time 
we set $k = 1$. Let  $\bar{D}_1$ denote the proposition that $C_1$ is not dirty (it belongs to the set of facts) and set $s_\beta \leq \bar{D}_1$ where 
$C_1 \notin \beta$. Note that the situations in which $C_1$ is  not dirty satisfy this fact, 
for example $s_{\{C_1\}} \leq \bar{D}_1$.  Denote by $\bar{q}$ the first round of answers that includes child 
one's lying and the others' ``No!" replies. The appearance of this action to $C_1$ is the identity since he knows 
that he is lying $f_{C_1} (\bar{q}) = \bar{q}$, whereas other children who do not know that $C_1$ is lying think 
that the action $q$ in classical muddy children (truthful public refutation) is happening, that is for $1 < i \leq n$ 
we have $f_{C_i}(\bar{q}) = q$. The kernel of $\bar{q}$ is the downset of the proposition in which $C_1$ knows 
he is not dirty and others know that  they  are dirty i.e.
\[
Ker(\bar{q}) = \downarrow (\Box_{C_1} \bar{D}_1 \vee \bigvee_{i = 2}^{n} \Box_{C_i} D_i)\,.
\]

\begin{proposition} 
After the first child's lying and the others' negative answers in the first round, 
every clean child $j$ (with $1< j\leq k$) thinks (wrongly) that he is dirty   
i.e. 
\[
s_{\{C_1\}} \leq [q_0 \,\bullet \bar{q}] \Box_{C_j}D_j\,.   
\]
\end{proposition} 

\noindent 
{\bf Proof.}  
We proceed in the same way as above. By moving the dynamic and epistemic modalities to the left and applying the update inequality eq.(\ref{eq:update}) we obtain
\[
f_{C_j}(s_{\{C_1\}}) \cdot f_{C_j} (q_0) \cdot  f_{C_j}(\bar{q}) \leq D_j\,.
\]
By replacing the $f_{C_j}$'s with their values we get 
\[
(s_{\{C_1\}} \vee s_{\{C_1, C_j\}}) \cdot  q_0 \cdot  q \leq D_j
\]
and by distributivity  we have to show the following two cases (the same as in the classical version above) 
\[
s_{\{C_1\}} \cdot  q_0 \cdot  q \leq D_j\qquad\quad and\qquad\quad s_{\{C_1, C_j\}} \cdot  q_0 \cdot  q \leq D_j \,.
\]
The second case is trivial for the same reasons as classical muddy children. For the first case we use eq.(\ref{eq:kernel}) proved by induction above and get $s_{\{C_1\}} \cdot  q_0 \in Ker(q)$ and hence $\bot = s_{\{C_1\}} \cdot  q_0 \cdot  q \leq D_j\,.$
\endproof\newline 

\paragraph{Secret Communication.}
As another example, consider the original $n$ and $k$ version but in which, just before the $k-1$'th round,
all but one of the dirty children (say, all except $C_1$), ``cheat" by secretly telling each other that they are in   
fact dirty. We denote this action as $\pi$. In the $k-1$'th round, all these dirty cheating children will announce that they know they are dirty 
(or equivalently refute that they do not know that they are dirty) where as $C_1$ answers as usual. We denote this mixed round of answers by $q'$. For the encoding of these actions in epistemic systems, that is their appearance and kernels refer to~\cite{BaltagSadr}. Now following the same line as in the proofs above, we can show that in the $k$'th round the only non-cheating child $C_1$ will wrongly conclude that he is clean i.e.
\[
s_{C_1,\dots,C_k} \leq [q_0 (\bullet \, q)^{k-2} \bullet \pi \bullet q'] \Box_1\bar{D_1}\,.
\]
The proof is done similar to the above cases and is presented in detail in~\cite{BaltagSadr}. 

\paragraph{A cryptographic attack.}  
This cryptographic attack is a somewhat simplified version of the man in the middle (MITM) attack which is a primary defect of public key-based systems. Two agents $A$ and $B$ share a secret key so 
that they can send each other encrypted messages over some 
communication channel. The channel is not secure: some outsider $C$ 
may intercept the messages or prevent them from being delivered 
(although he cannot read them because he does not have the 
key). Suppose the encryption method is publicly known but the key is 
secret. It is also known that $A$ is the only one who knows an important 
secret for example if some fact $P$ holds or not. Suppose now that $A$ 
sends an encrypted message to $B$ communicating the secret. $B$ gets 
the message and he is convinced that it must be authentic. Now both $A$ and $B$ 
are convinced that they share the secret and that $C$ doesn't. However 
suppose that $C$ notices two features of the specific encryption 
method: first that the shape of the encrypted message can show whether 
it contains a secret or it is just junk, second that without knowing 
the key or the content of the message he can modify the encrypted 
message to its opposite i.e. if it originally said $P$ holds, it will 
now say that $P$ does not hold. The outsider $C$ will then secretly 
intercept the message, change it appropriately and send it to $B$ 
without knowing the secret. Now 
$A$ and $B$ mistakenly believe that they share the secret, while in 
fact $B$ got the wrong secret instead and  $C$ has succeeded to manipulate 
their beliefs. 
 
We can encode this situation in an epistemic system. 
The agents include $\{A,B,C\}$ and we call  the message in which $P$ holds  $P$ and the 
one in which it does not hold  $\bar{P}$, these are inconsistent facts so we have  $P,\bar{P}\in Stab(Q)$ and $P \wedge \bar{P} = \bot, 
\ P \vee \bar{P} = \top$. Let $s,t\in M$ 
satisfy  $s \leq P$ and $t \leq \bar{P}$. The only agent that 
knows if $P$ holds or not is $A$ thus  $f_A(s) = s$ and similarly 
$f_A(t) = t$. On the other hand $B$ and $C$ do not know this so $f_B(s) = f_C(s) = f_B(t) = f_C(t) = 
s\vee t$.  
 The epistemic actions 
that correspond to the cryptographic attack are 
the following: $\alpha$ in which the message $P$  is intercepted, modified and sent 
to $B$, $\beta$ in which the message $\bar{P}$ is intercepted, 
modified and sent to $B$, $\alpha'$ in which $A$ sends the message $P$ 
 to $B$, $\beta'$ in which $A$ sends the message $\bar{P}$ to $B$, and finally $\gamma$  which corresponds to sending a 
junk message. Thus $\{\alpha, \beta, \alpha', \beta', \gamma\} 
\subseteq Q$. In actions $\alpha$ and 
$\beta$ agent $C$ is uncertain about which message $P$ or $\bar{P}$ has been 
sent so $f_C(\alpha) = f_C(\beta) = \alpha\vee\beta$. On the other hand, 
agent $A$ is sure that he has sent a message (either that $P$ holds or 
that it doesn't) to $B$ and that $B$ has 
received exactly the same secret i.e. 
$f_A(\alpha) = \alpha'$ and $f_A(\beta) = \beta'$. However if 
$P$ has been sent, $B$ has received $\bar{P}$ so $f_B(\alpha) = \beta'$ 
and the other way around $f_B(\beta) = \alpha'$. Further $f_A(\alpha') 
= f_B(\alpha') = \alpha'$ and $f_A(\beta') = f_B(\beta') = 
\beta'$ and 
$f_C(\alpha') = f_C(\beta') = \alpha'\vee\beta'\vee\gamma$. $C$ also 
considers it possible that  only a junk message has been sent and that is why 
he sees $\gamma$ while in $\alpha'$ and $\beta'$. If a junk message 
has been sent, $A$ and $B$ are sure about it $f_A(\gamma) = 
f_B(\gamma) = \gamma$ while $C$ is unsure if it was a junk message 
or $P$ or $\bar{P}$, thus $f_C(\gamma) = 
\alpha' \vee \beta' \vee \gamma$. The kernel of each action comprises the states 
which they cannot be applied to i.e. $Ker(\alpha)= Ker(\alpha') = \ \downarrow \!\bar{P}$ and $Ker(\beta) = 
Ker(\beta') = \ \downarrow \!P$. 
 
The epistemic action $\alpha \vee \beta$ expresses the action of  
communicating the secret $P$ or $\bar{P}$ in the above scenario. Now let us update  
the state $s$ with the epistemic action $\alpha \vee  
\beta$ and  
show that after update, if $P$ holds, then $A$ knows that $B$ knows that  
$P$ holds, that is
\[
s \cdot  (\alpha \vee \beta) \leq \Box_A \Box_B P
\]
 Since this is equal to  $(s \cdot  \alpha)  
\vee (s \cdot \beta) \leq \Box_A \Box_B P$ and $s\leq P \in Ker(\beta)$ we  
get $s \cdot  \beta = \bot$, so it suffices to show that $s \cdot   
\alpha \leq \Box_A \Box_B P$. By adjunction $f_B(f_A(s \cdot   
\alpha)) \leq P$. By eq.(\ref{eq:update}) we get 
$f_A(s \cdot  \alpha) \leq f_A(s) \cdot  f_A(\alpha)$, and 
order preservation of $f_B$ will give us
\[f_B(f_A(s \cdot  \alpha)) \leq f_B(f_A(s) \cdot  f_A(\alpha)) \leq f_B(f_A(s)) \cdot  f_B(f_A(\alpha)).
\]
Now it suffices to show $f_B(f_A(s)) \cdot  f_B(f_A(\alpha)) \leq P$.
We do that by replacing $f_A$ with its values and show $f_B(s) \cdot  f_B(\alpha') \leq P$, do the same for $f_B$ and get $(s \vee t) \cdot  \alpha' \leq P$,  
hence   
$(s \cdot  \alpha') \vee (t \cdot  \alpha') \leq P$ which is equal to $(s \cdot  \alpha') \leq P$ since  
$t \leq \bar{P}\in Ker(\alpha')$.  By the assumption  $s \leq P$ we obtain $s \cdot  \alpha' \leq P \cdot  
\alpha'$ which leads to $s \cdot  \alpha' \leq P$ because $P$ is a fact.  

\paragraph{{\bf BMS} Models as Epistemic Systems.} The Kripke semantics for dynamic epistemic logic as introduced in~\cite{BaltagMossSolecki, BaltagMoss} are examples of epistemic systems by the following theorem:
\begin{theorem}\label{thm:represe}
\em  Models of  {\bf BMS} are epistemic systems $(M,Q,\{f_A\}_{A \in {\cal A}})$ with the following properties
\begin{enumerate}
\item Both $M$ and $Q$ are completely distributive atomistic Boolean algebras.
\item If $m$ is an atom of $M$ and $q$ is an atom of $Q$ then $m \cdot  q$ is either $\bot$ or an atom. 
\item If $q,q'$ are atoms of $Q$ then $q \bullet q'$ is an atom. 
\item If $m \cdot  q = \bot$ then either $m = \bot$ \ or \ $q = \bot$.
\end{enumerate}
\end{theorem}  
The proof goes by constructing an epistemic system given a model of {\bf BMS} and is presented in detail  in \cite{SadrThesis} --- and is based on ideas introduced in \cite{BCS}.  Key is the observation that each relation $R \subseteq S \times S$  gives rise to a sup-map $f_R: {\cal P}(S) \to {\cal P}(S)$ i.e.~we  lift the accessibility relations $\to_A$ of  the Kripke semantics of {\bf BMS} to appearance maps $f_A$.  A model of {\bf BMS} consists of two Kripke structures, one for the states as usual $(S, \to_A, \mu)$ and one for the deterministic actions  $(\Sigma, \to_A, pre)$ where $pre: \Sigma \to {\cal P}(S)$ assigns to each actions a precondition. The state model acts on the action model resulting in an updated state model, via a partial cartesian product, that is the \emph{epistemic update}. 
Action models act on themselves via a sequential composition operation. 
In order to construct an epistemic system, we close the set of `states' (= deterministic actions) of the action model under sequential composition  and close the set of states of the state model under  update. The closure of the states yields the atoms of module and the closure of deterministic actions yields the atoms of quantale, and we  get a Boolean epistemic system by taking their  powersets $({\cal P}(S),{\cal P}(\Sigma), \{f_A\}_{A \in \cal A})$.  Operations of this epistemic system are constructions of {\bf BMS}, e.g.~epistemic update and  sequential composition   extended pointwisely to subsets of states and actions. The  epistemic and dynamic modalities arise, as before, as adjoints to the lifted appearance and update maps, but moreover and because of the boolean complementation we get  a de Morgan dual for each of these modalities, in particular the de Morgan dual of the epistemic modality $(\Box_A (-)^c)^c$ stands for the $\Diamond$-modality of standard epistemic logic.

\section{The sequent calculus of epistemic systems}  
We have two different sequent systems,  a $Q$-system and an $M$-system, both of them are intuitionistic in the sense that they have only one formula on the right hand side of the turnstile. The $Q$-system corresponds to the quantale part and  the $M$-system corresponds to the module of a \emph{distributive   epistemic system}. By this we mean   an epistemic system with a distributive module. 

\paragraph{The $Q$-system.} The formulas of the $Q$-system, denoted as $L_Q$, are generated by the following syntax:
\[
q ::= \top \mid \bot \mid \sigma \mid 1 \mid q \bullet q \mid q\,/\,q \mid q\setminus q  \mid q \vee q  \mid q \wedge q \mid f_A^Q(q) \mid \Box_A^Q\,q
\]
where $A$ is in the set ${\cal A}$ of agents,  and $\sigma$ is in a set $V_Q$ of atomic  
action variables. The sequents of the $Q$-system are called $Q$-sequents and are denoted as
 \[
\Gamma \vdash_Q q
\]
where $\Gamma$ is a sequence of actions and agents, that is $\Gamma \in (L_Q \cup {\cal A})^*$, and $q$ is a single action, that is $q \in L_Q$.  To assign meaning to the sequents of the $Q$-system, we introduce
\[
- \odot_Q -: L_Q \times (L_Q \cup {\cal A}) \to L_Q
\]
by putting 
\[
q \odot_Q q' := q \bullet q'\qquad\qquad\qquad\qquad q \odot_Q A := f_A^Q (q)\,.
\]
For $\Gamma = (\gamma_1, \cdots, \gamma_n) \in   (L_Q \cup {\cal A})^*$
we take the convention
\[
\bigodot_Q \Gamma := ((((1 \odot_Q \gamma_1) \odot_Q  
\gamma_2)\odot_Q \gamma_3)\cdots)\odot_Q \gamma_n\,.
\]
As an example, the sequence $\Gamma = (q, A, q')$ corresponds to 
\[
\bigodot_Q \Gamma = ((1 \odot_Q q) \odot_Q A) \odot_Q q' = f_A^Q(1 \bullet q) \bullet q'  = f_A^Q(q) \bullet q'\,.
\]
Adding the multiplicative unit to the beginning of the sequence will allow us to avoid non-well defined $\odot_Q$-expressions for sequences such as $\Gamma = A$, which will now mean $\bigodot_Q A = f_A^Q(1)$.  Indeed, the operation 
$\bigodot_Q$ constitutes our semantic interpretation of the comma for $Q$-sequents.
For simplicity we denote the semantics of a formula by the formula itself i.e.~rather than $\semantics{\bigodot_Q \Gamma}$ and $\semantics q$ we write $\bigodot_Q \Gamma$ and $q$. 
The empty sequence on the left hand side stands for 1 and we do not allow for the empty sequence on the right hand side. We define a \emph{satisfaction} relation $\models_Q$  on the $Q$-system as follows:
\[
\Gamma \models_Q q'\ \ \Longleftrightarrow\ \ \bigodot_Q \Gamma \leq q'\,.
\]
We say that a sequent  $\Gamma \vdash_Q q'$ is \emph{valid} if and only if $\Gamma \models_Q q'$.
In this way we identify any $Q$-sequence $\Gamma$ with a $Q$-formula and  its corresponding element of the quantale. 

\smallskip
Ordered monoids have first  been used by Lambek to model Lambek-calculus. Yetter showed  that quantales are  models of  non-commutative Linear Logic \cite{Yetter}. The extension of these systems to epistemic modalities and quantales with operators is new.  So the operational and unit rules for the $Q$-system are the rules for Non-Commutative Intuitionistic Linear Logic, extended with an agent context. In order to see the connection with Linear Logic  note that our multiplication $\bullet$ is the tensor of Linear Logic, our disjunction is the Linear Logic  sum, the conjunction is $\&$, and our left and right residuals are $\circ -$ and $-\circ$. In a table:
\begin{center}
\begin{tabular}{|c|c|}
\hline
Q-system & Linear Logic\\
\hline
\hline
1 & 1\\
\hline
$\top$ & $\top$\\
\hline
$\bot$ & 0\\
\hline
$\bullet$ & $\otimes$\\
\hline
$/$ & $\circ -$\\
\hline
$\setminus$ & $-\circ$\\
\hline
$\vee$ & $\oplus$\\
\hline
$\wedge$ & $\&$\\
\hline

\end{tabular}
\end{center}
The axiom and unit rules of the $Q$-system are:
\begin{center}
\fbox{$
\begin{array}{ccc}
\infer[Id]{q \vdash_Q q}{}\hspace{1cm}
&
\infer[1L]{\Gamma, 1, \Gamma' \vdash_Q q}{\Gamma, \Gamma' \vdash_Q q}\hspace{1cm}
&
\infer[1R]{\vdash_Q 1}{} 
\\ \ \\
\infer[\bot L]{\Gamma, \bot, \Gamma' \vdash_Q q}{} \hspace{2cm}
&
\infer[\top R]{\Gamma \vdash_Q \top}{}
\end{array}
$}
\end{center}
The operational rules of the $Q$-system including those on epistemic modalities are:
\begin{center}
\fbox{$
\begin{array}{cc}
\infer[f_A^Q R]{\Gamma, A \ \vdash_Q \ f_A^Q(q)}{\Gamma \ \vdash_Q \ q}
&
\infer[f_A^Q L]{f_A^Q(q'), \Gamma \ \vdash_Q \ q}{q', A,   \Gamma \ \vdash_Q \ q}\\ \ \\
\infer[\Box_A^Q R]{\Gamma \vdash_Q \Box_A^Q\, q}{\Gamma, A \vdash_Q q}
&
\infer[\Box_A^Q L]{\Box_A^Q\, q', A, \Gamma \vdash_Q q}{q', \Gamma \vdash_Q q}
\\ \ \\
\infer[\bullet L]{\Gamma, q_1 \bullet q_2, \Gamma' \vdash_Q q}{\Gamma, q_1, q_2, \Gamma' \vdash_Q q} 
& \infer[\bullet R]{\Gamma_Q, \Gamma'_Q, \Gamma_A \vdash_Q q_1 \bullet q_2}{\Gamma_Q, \Gamma_A \vdash_Q q_1 \quad \Gamma'_Q, \Gamma_A \vdash_Q q_2}
\\ \ \\
\infer[\vee L]{\Gamma, q_1 \vee q_2, \Gamma' \vdash_Q q}{\Gamma, q_1, \Gamma' \vdash_Q q \qquad \Gamma, q_2, \Gamma' \vdash_Q q} & \infer[\vee R1]{\Gamma \vdash_Q q_1 \vee q_2}{\Gamma \vdash_Q q_1}\hspace{0.6cm}\infer[\vee R2]{\Gamma \vdash_Q q_1 \vee q_2}{\Gamma \vdash_Q q_2}\\ \ \\
\infer[\wedge L1]{\Gamma, q_1 \wedge q_2, \Gamma' \vdash_Q q}{\Gamma, q_1, \Gamma' \vdash_Q q}\qquad
 \infer[\wedge L2]{\Gamma, q_1 \wedge q_2, \Gamma' \vdash_Q q}{\Gamma, q_2, \Gamma' \vdash_Q q}&
\infer[\wedge R]{\Gamma \vdash_Q q_1 \wedge q_2}{\Gamma \vdash_Q q_1 \qquad \Gamma \vdash_Q q_2}\\ \ \\
\infer[/ L]{ q_1\,/\,q_2, \Gamma_Q  \vdash_Q q}{\Gamma_Q \vdash_Q q_2 \quad q_1 \vdash_Q q}  &
\infer[/R]{\Gamma \vdash_Q q_1\,/\,q_2}{\Gamma, q_2 \vdash_Q q_1}
\\ \ \\
\infer[\setminus L]{ \Gamma, q_1\setminus q_2 \vdash_Q q}{\Gamma \vdash_Q q_1 \quad  q_2 \vdash_Q q}&
\infer[\setminus R]{\Gamma_Q \vdash_Q q_1\setminus q_2}{q_1, \Gamma_Q \vdash_Q q_2}

\end{array}    
$}
\end{center}
As in non-commutative Linear Logic we have no weakening, contraction and exchange rules for actions. Our structural rules consist of a restricted version of the  usual cut rule and a rule to encode the relation between appearance maps and the unit of composition --- eq.(\ref{eq:skip}) in the algebra. These two rules are: 
\begin{center}
\fbox{$
\begin{array}{cc}
\infer[Qcut]{\Gamma', \Gamma'' \vdash_Q q'}{\Gamma' \vdash_Q q \qquad  q, \Gamma'' \vdash_Q q'}\hspace{2cm} & \infer[Agent]{1 \vdash_Q q}{A \vdash_Q q}
\end{array}
$}
\end{center}

\paragraph{The $M$-system.} The formulas of the $M$-system, denoted as $L_M$, are generated by the following syntax:
\[
m ::= \bot \mid \top \mid p \mid s \mid m \wedge m \mid m \vee m  \mid [q]m \mid m \cdot  q \mid  
\Box_A^M\, m \mid f_A^M(m)
\]
where $A$ is in the set ${\cal A}$ of agents, $p$ is in the set $\Phi$ of facts, and $s$ is in a set  
$V_M$  of atomic propositional variables. The sequents of the $M$-system are called $M$-sequents and are denoted as 
\[
\Gamma \vdash_M m
\]
where $\Gamma$ is a sequence of propositions, actions, and agents, that is   $\Gamma \in (L_M \cup L_Q \cup {\cal A})^*$ and  $m$ is a single proposition, that is $m \in L_M$. 
To assign meaning to the sequents of an $M$-sequent we consider   
\[  
- \odot_M -: L_M \times (L_M \cup L_Q \cup {\cal A})\rightarrow 
L_M 
\]  
now by putting 
\[
m \odot_M A := f_A^M(m)\qquad\qquad m \odot_M q := m  \cdot  q\qquad\qquad m \odot_M m' := m \wedge m'\,.
\] 
For $\Gamma = (\gamma_1, \cdots, \gamma_n) \in  
(L_M \cup L_Q \cup {\cal A})^*$ we take the convention 
\[
\bigodot_M \Gamma := ((((\top \odot_M \gamma_1) \odot_M  
\gamma_2)\odot_M \gamma_3)\cdots)\odot_M \gamma_n\,.
\]
As an example, the sequence $\Gamma = (m, A, q, B, m')$ corresponds to 
\[
\bigodot_M \Gamma = ((((\top \odot_M m) \odot_M A) \odot_M q) \odot_M B) \odot_M m'= f_B(f_A(\top \wedge m) \cdot  q) \wedge m' = f_B(f_A(m) \cdot  q) \wedge m'\,.
\] 
Sequences  of only one agent $\Gamma = A$ will mean $\bigodot_M \Gamma = f_A^M(\top)$ and sequences of only one action $\Gamma = q$ will mean $\bigodot_M \Gamma = \top\cdot q$. The empty sequence on the left hand side stands for $\top$, here we also allow for an empty right hand side, which stands for $\bot$. As before,  we denote the semantics of a formula by the formula itself.
We define 
\[
\Gamma \models_M m'\ \ \Longleftrightarrow\ \ \bigodot_M \Gamma \leq m'\,.
\]
and  say that a sequent $\Gamma \vdash_M m'$ is  \emph{valid} if and only if \ $\Gamma \models_M m'$.  In this way we identify any $M$-sequence $\Gamma$ with an $M$-formula and its  corresponding element of the module.    

The rules of the $M$-system correspond to a  distributive lattice logic extended with an agent context for our epistemic modalities. The axiom and unit rules of the $M$-system are\footnote{The $\bot R$ rule follows from the $weakR$ rule and thus can be dropped.}:
\begin{center}
\fbox{$
\begin{array}{cccc}
\infer[Id]{m \vdash_M m}{}\hspace{0.5cm}&
\infer[\bot L]{\bot \vdash_M m}{} \hspace{0.5cm}&
\infer[\bot R]{\Gamma \vdash_M \bot}{\Gamma \vdash_M \quad}\hspace{0.5cm} & \infer[\top R]{\Gamma \vdash_M \top}{}
\end{array}
$}
\end{center}

\medskip \noindent
The operational rules of the $M$-system for the lattice operations and modalities are: 
\begin{center}
\fbox{$
\begin{array}{cc}
\infer[f_A^M R]{\Gamma, A \ \vdash_M \ f_A^M(m)}{\Gamma \ \vdash_M \ m}
&
\infer[f_A^M L]{f_A^M(m), \Gamma \ \vdash_M \ m'}{m, A,  \Gamma \ \vdash_M \ m'}
\\ \ \\
\infer[\Box_A^M R]{\Gamma \vdash_M \Box_A^M\,m}{\Gamma, A \vdash_M  m}
&
\infer[\Box_A^M L]{\Box_A^M (m), A, \Gamma \vdash_M m'}{m, \Gamma \vdash_M m'}
\\ \ \\
\hspace{-2ex}\infer[\wedge R]{\Gamma \vdash_M m_1 \wedge m_2}{\Gamma \vdash_M m_1 \quad \Gamma \vdash_M m_2} & \hspace{-2ex}\infer[\wedge L1]{\Gamma, m_1 \wedge m_2, \Gamma' \vdash_M m}{\Gamma, m_1, \Gamma' \vdash_M m} \hspace{0.4cm}
\infer[\wedge L2]{\Gamma, m_1 \wedge m_2, \Gamma' \vdash_M m}{\Gamma, m_2, \Gamma' \vdash_M m}
\\ \ \\
\infer[\vee L]{\Gamma, m_1 \vee m_2, \Gamma'}{\Gamma, m_1, \Gamma' \vdash_M m \quad \Gamma, m_2, \Gamma' \vdash_M m}&
\infer[\vee R1]{\Gamma \vdash_M m_1 \vee m_2}{\Gamma \vdash_M m_1}
\hspace{0.7cm}
\infer[\vee R2]{\Gamma \vdash_M m_1 \vee m_2}{\Gamma \vdash_M m_2}
\end{array}
$}
\end{center}
As structural rules, we have  propositional weakening, contraction, and exchange,  a restricted version of the usual cut  rule, and a rule encoding stability of facts under update:
\begin{center}
\fbox{$
\begin{array}{ccc}
\infer[contr]{\Gamma, m', \Gamma' \vdash_M m}{\Gamma, m', m', \Gamma' \vdash_M m}&
\infer[exch]{\Gamma, m', m'', \Gamma' \vdash_M m}{\Gamma, m'', m', \Gamma' \vdash_M m}&
\infer[fact]{\Gamma, q \vdash_M p}{\Gamma \vdash_M p}
\\ \ \\
\infer[Mcut]{\Gamma', \Gamma'' \vdash_M m}{\Gamma' \vdash_M m' \qquad  m', \Gamma'' \vdash_M m}
&
\infer[weakL]{\Gamma, m', \Gamma' \vdash_M m}{\Gamma, \Gamma' \vdash_M m}
&
\infer[weakR]{\Gamma \vdash_M m}{\Gamma \vdash_M}
\end{array}
$}
\end{center}

\paragraph{The $MQ$-system.}
Since the core of our approach is the action of the quantale on the module, we also have \em mixed rules \em for epistemic update and dynamic modality consisting of both $M$- and $Q$-sequents: 
\begin{center}
\fbox{$
\begin{array}{cc}
\infer[\cdot  L]{m' \cdot q, \Gamma \vdash_M m}{m',  q,  \Gamma \vdash_M m}\hspace{1cm}
&
\infer[\cdot  R]{\Gamma, \Gamma_Q, \Gamma_A \vdash_M m\cdot q}{\Gamma , \Gamma_A \vdash_M m & \Gamma_Q, \Gamma_A \vdash_Q q}
\\ \ \\
\infer[Dy L]{ [q]m', \Gamma_Q \vdash_M m}{m' \vdash_M m \quad  \Gamma_Q \vdash_Q q}\hspace{1cm}
&
\infer[Dy R]{\Gamma \vdash_M [q]m}{\Gamma, q \vdash_M m}
\end{array}
$}
\end{center}
The action of the quantale on the module  preserves the unit of multiplication, is  disjunction preserving in both arguments, and satisfies an associativity condition with regard to composition of actions. In order to prove the same properties  for  epistemic update (and dual ones for dynamic modality) in the $M$-system,  we should be able to work with the quantale operations in $M$-sequents.  So we have the following rules that include for example update with unit, composition, and choice of actions:
\begin{center}
\fbox{$
\begin{array}{cc}
\infer[1ML] {\Gamma, 1, \Gamma' \vdash_M m}{\Gamma, \Gamma' \vdash_M m}\hspace{2cm}
&
\infer[\bullet ML]{\Gamma, q_1 \bullet q_2, \Gamma' \vdash_M m}{\Gamma, q_1, q_2, \Gamma' \vdash_M m}
\\ \ \\
 \infer[\vee ML]{\Gamma, q_1 \vee q_2 \vdash_M m}{\Gamma, q_1 \vdash_M m \qquad \Gamma, q_2 \vdash_M m}
\\ \ \\
\infer[/ ML]{\Gamma, q_1\,/\,q_2, \Gamma_Q \vdash_M m}{\Gamma_Q \vdash_Q q_2 \quad \Gamma, q_1 \vdash_M m}\hspace{2cm} &
\infer[\setminus ML]{\Gamma, \Gamma_Q, q_1\setminus q_2 \vdash_M m}{\Gamma_Q \vdash_Q q_1 \quad \Gamma, q_2 \vdash_M m}\\ \ \\
\infer[\wedge ML1]{\Gamma, q_1 \wedge q_2, \Gamma' \vdash_M m}{\Gamma, q_1, \Gamma' \vdash_M m}\hspace{2cm}&
\infer[\wedge ML2]{\Gamma, q_1 \wedge q_2, \Gamma' \vdash_M m}{\Gamma, q_2, \Gamma' \vdash_M m}
\end{array}
$}
\end{center}

\paragraph{Note on the cut rules.} 
We have two cut rules in our system: a $QCut$ for the $Q$-system and an $MCut$ for the $M$-system. 
Although  the $Q$-system subsumes  a quantale logic for which cut is eliminable, for example non-commutative Linear Logic  or Lambek-calculus, the $QCut$  of our system does not inherit  this property.  We believe this is partly because of the modal part of the logic, that is the quantale endomorphisms   and their  interaction with the non-commutative sequential composition in eq.~(\ref{eq:mult}). This equation is encoded in the  $\bullet R$ rule, which  needs context splitting for actions and context sharing for agents.  
Similarly, the $M$-cut is not eliminable,    partly due to the update inequality eq.~(\ref{eq:update}) and its corresponding $\cdot  R$ rule.   However (as noticed by  one of our referees), in both of these systems the identity rules can be reduced to atomics, which is a  sign of well-definedness of our system. Studying these proof theoretic issues constitutes future work.

\paragraph{Note on intuitive reading of sequents.} To provide the reader with a 
way to read our sequents in natural language, we capture the {\it intuitive meaning} of a sequent in the following
inductive manner:
\begin{itemize}
\item $\vdash_M m$ means ``proposition $m$ holds in all contexts''
\item $\vdash_Q q$  means ``action $q$  does not necessarily have an effect on propositions''
\item $\Gamma, A, \Gamma' \vdash_M  m$ means ``in context $\Gamma$ agent $A$ knows/believes, that $\Gamma'\vdash_M m$ holds'' --- this
captures features of $A$'s own reasoning: $\Gamma'\vdash_M
m$ is accepted by $A$ in context $\Gamma$ as a valid argument.
\item $\Gamma, q, \Gamma' \vdash_M m$ means ``after action $q$ happens on context $\Gamma$, the sequent  $\Gamma'\vdash_M m$ will hold''
\item $m, \Gamma \vdash_M m$ means ``in context $m$
(in any situation in which $m$ is true), the sequent
 $\Gamma\vdash_M m$ holds''
\item $\Gamma, A, \Gamma' \vdash_Q q$ means ``in context $\Gamma$ agent $A$ knows/believes,  that $\Gamma' \vdash_Q q$ holds''
\item $q, \Gamma \vdash_Q q$ means ``after doing action $q$, the sequent $\Gamma \vdash_Q q$ holds''
\end{itemize}
Observe that the left-to-right order of this  intuitive reading is the opposite  of the right-to-left application order of $\odot$ or comma. This is because the reading
involves the (intuitive) notions of \em knowledge \em and \em weakest precondition\em, which are adjoints of 
the $f_A$ and $\cdot$ operations; thus, the intuitive reading can be obtained by taking the adjoints
(which live on the right-side of turnstile) of the formulas on the left-hand side of a sequent. 
For instance, the sequent $m, A, B \vdash_M m'$ after applying commas on the left would mean  $f_B^M(f_A^M(m)) \leq m'$, and after applying the adjoints would correspond to $m \leq \Box_A^M  \Box_B^M m'$. This has now the exact shape of its intuitive meaning which is  ``in context $m$ agent $A$ believes that agent $B$ believes that $m'$''. Examples such as $\Gamma, m \vdash_M m'$ make more sense when $M$ is a Heyting algebra. For instance, the sequent $m, A, q, B, m' \vdash_M m''$ can be read
as:
``in context $m$, agent $A$ believes that after action $q$ agent $B$
will
believe that, in context $m'$, proposition $m''$ must hold''.
This reading shows that, as already mentioned in the introduction,  our sequent calculus expresses two forms of resource sensitivity. 
One is the  use-only-once form of linear logic~\cite{Girard} that comes from the quantale structure on epistemic actions, which we called \em dynamic resources\em. 
The other form deals with \em epistemic resources\em: the resources available to each agent that enable him to  
reason in a certain way (i.e. to infer a conclusion from some assumptions). These resources are encoded in the way 
the context appears to the agent in sequents, for instance $\Gamma$ in the sequent $\Gamma, A, \Gamma' \vdash_M m$ 
is the context,  and hence $f_A^M(\Gamma)$ is the resource that  enables agent $A$ to do the $\Gamma' \vdash_M m$ 
reasoning. Note that $\Gamma' \vdash_M m$ might not be a valid sequent in the context $\Gamma$, but it is valid 
in the context given by  $\Gamma$'s  appearance  to agent $A$. 

\medskip
\begin{theorem}[Soundness] \em The rules presented in this section are sound  with respect to the 
algebraic semantics in terms of distributive epistemic systems.
\end{theorem}

\paragraph{Proof.} For the soundness of the rules we have to show that derivable sequents  of the $Q$ and $M$-systems are valid in a distributive epistemic system, that is $Q$-rules are valid in the quantale part and $M$-rules (including the mixed rules) are valid in the module part. In other words, we have to prove the following for the $Q$-system 
\[
\mbox{if} \qquad \Gamma \vdash_Q q \qquad \mbox{then} \qquad \Gamma \models_Q q
\]
and  a similar one for the $M$-system.  The proof is done  by induction on the $\odot$ operation and 
applying the algebraic definitions and properties of the connectives to show  that the rules of each system preserve validity of sequents. 
The full proof  can be found in~\cite{SadrThesis},   here we provide the reader with   the proofs for  the rules that show the crucial features of our system:  sequential composition and appearance and knowledge of actions in the quantale for the $Q$-system, and  epistemic update and dynamic modality for the $M$-system.

\smallskip
\noindent{\bf i.~Soundness of sequential composition.} The rules for sequential composition are
\[
\infer[\bullet L]{\Gamma, q_1 \bullet q_2, \Gamma' \vdash_Q q}{\Gamma, q_1, q_2, \Gamma' \vdash_Q q} \hspace{3cm} \infer[\bullet R]{\Gamma_Q, \Gamma'_Q, \Gamma_A \vdash_Q q_1 \bullet q_2}{\Gamma_Q, \Gamma_A \vdash_Q q_1 \quad \Gamma'_Q, \Gamma_A \vdash_Q q_2}
\]
To prove soundness, we have to show that if the sequent on the top line is valid, so is the sequent on the bottom line. Using the definition of validity, we have to show the following satisfaction statement for the left rule:
\[
\mbox{If} \quad  \Gamma, q_1, q_2, \Gamma' \models_Q q \quad \mbox{then} \quad  \Gamma, q_1 \bullet q_2, \Gamma' \models_Q q
\]
By the definition of satisfaction in terms of $\odot_Q$,  we have to show the following 
\[
\mbox{If} \quad \bigodot_Q (\Gamma, q_1, q_2, \Gamma')\leq q \quad \mbox{then} \quad \bigodot_Q(\Gamma, q_1 \bullet q_2, \Gamma') \leq q\,,
\]
This is true since the application of  $\bigodot_Q$ to the left hand side sequences of the top and bottom sequents yields equal quantale elements, that is 
 \[
 \bigodot_Q (\Gamma, q_1, q_2, \Gamma') =  \bigodot_Q\!\Gamma \bullet q_1 \bullet q_2 \bullet \bigodot_Q\!\Gamma' = \bigodot_Q(\Gamma, q_1\bullet q_2, \Gamma')
 \]
For the right rule we proceed similarly and show the following satisfaction statement 
\[
\mbox{If} \quad \Gamma_Q, \Gamma_A \models_Q q_1 \quad \mbox{and} \quad \Gamma'_Q, \Gamma_A \models_Q q_2 \quad \mbox{then} \quad \Gamma_Q, \Gamma'_Q, \Gamma_A \models_Q q_1 \bullet q_2
\]
which is by definition equivalent  to the following  $\bigodot_Q$ statement
\[
\mbox{If} \quad \bigodot_Q(\Gamma_Q, \Gamma_A) \leq  q_1 \quad \mbox{and} \quad \bigodot_Q (\Gamma'_Q, \Gamma_A) \leq q_2 \quad \mbox{then} \quad \bigodot_Q(\Gamma_Q, \Gamma'_Q, \Gamma_A) \leq  q_1 \bullet q_2
\]
We first assume that we have only one agent in our agent context, that is $\Gamma_A = A$ and  we have to show the following 
\[
\mbox{If} \quad f_A^Q\Bigl(\bigodot_Q\Gamma\Bigr) \leq q_1 \quad \mbox{and} \quad f_A^Q\Bigl(\bigodot_Q\Gamma'_Q\Bigr) \leq q_2 \quad \mbox{then} \quad f_A^Q \Bigl(\bigodot_Q\Gamma_Q \bullet \bigodot_Q\Gamma'_Q\Bigr) \leq q_1\bullet q_2
\]
Assume that the precedent holds, 
 by order-preservation of the multiplication on the quantale we can multiply both sides of these inequalities and we get 
 \[
 f_A^Q\Bigl(\bigodot_Q\Gamma\Bigr) \bullet f_A^Q\Bigl(\bigodot_Q\Gamma'_Q\Bigr) \leq q_1 \bullet q_2\,, 
 \]
By the relation between appearance maps and multiplication on the quantale eq.(\ref{eq:mult}) we have 
\[
f_A^Q \Bigl(\bigodot_Q \Gamma \bullet \bigodot_Q \Gamma'_Q\Bigr) \leq f_A^Q\Bigl(\bigodot_Q\Gamma\Bigr) \bullet f_A^Q\Bigl(\bigodot_Q\Gamma'_Q\Bigr)\,,
\quad 
\mbox{hence}
\quad
f_A^Q \Bigl(\bigodot_Q \Gamma \bullet \bigodot_Q \Gamma'_Q\Bigr) \leq q_1 \bullet q_2\,.
\]
which is exactly what we wanted to prove, that is  the validity of the bottom line of the rule. If $\Gamma_A$ has more than one agent $\Gamma_A = A_1, \dots, A_n$ then we have to show that if  
\[ 
f_{A_1}^Q \Bigl(f_{A_2}^Q\Bigl(\dots f_{A_n}^Q\Bigl(\bigodot_Q\Gamma\Bigr)\Bigr)\Bigr) \leq q_1 \quad  \mbox{and} \quad  f_{A_1}^Q \Bigl(f_{A_2}^Q\Bigl(\dots f_{A_n}^Q\Bigl(\bigodot_Q\Gamma'_Q\Bigr)\Bigr)\Bigr) \leq q_2
\]
then
\[
f_{A_1}^Q \Bigl(f_{A_2}^Q\Bigl(\dots f_{A_n}^Q\Bigl(\bigodot_Q\Gamma_Q \bullet \bigodot_Q\Gamma'_Q\Bigr)\Bigr)\Bigr) \leq q_1\bullet q_2\,.
\]
The proof for this case is done similarly, except that after multiplying the two sides of the assumption by $\bullet$, we have to apply the inequality for $f_{A_i}^Q$ and the quantale multiplication $n$ times, that is once for each agent $A_i \in \Gamma_A$, starting from the innermost one $f_{A_n}^Q$ and ending with the outmost one $f_{A_1}^Q$. 

\smallskip
\noindent{\bf ii.~Soundness of appearance and knowledge of actions.}
The rules for the appearance map are
\[
\infer[f_A^Q L]{f_A^Q (q'), \Gamma \vdash_Q q}{q', A, \Gamma \vdash_Q q}  \hspace{3.5cm}
\infer[f_A^Q R]{\Gamma, A \vdash_Q f_A^Q (q)}{\Gamma \vdash_Q q}
\]
By using the satisfaction relation and definition of $\odot_Q$,  soundness of the left rule follows from definition of $\odot_Q$ between an agent and an  action. This is  because  $\bigodot_Q(q', A,\Gamma)$ is equal to $f_A^Q (q') \bullet \bigodot_Q\!\Gamma$, for which by the top line we  have $f_A^Q (q') \bullet \bigodot_Q\!\Gamma \leq q$. 
The right rule follows by the order preservation of $f_A^Q$, that is if $\bigodot_Q\!\Gamma \leq q$ then we have $f_A^Q (\bigodot_Q\!\Gamma) \leq f_A^Q (q)$, which is the meaning of the bottom line. 

\smallskip\noindent
The rules for knowledge on the quantale are:
\[
\infer[\Box_A L]{\Box_A^Q \,q', A, \Gamma \vdash_Q q}{q', \Gamma \vdash_Q q} \hspace{3.5cm}
\infer[\Box_A R]{\Gamma \vdash_Q \Box_A^Q \, q}{\Gamma, A \vdash_Q q}
\]
For the left rule assume $q' \bullet \bigodot_Q\!\Gamma \leq q$, and we have to show $f_A^Q (\Box_A^Q \, q') \bullet \bigodot_Q\!\Gamma \leq q$. By   composition  of adjoints on the $f_A^Q$ and $\Box_A^Q$,  we have  $f_A^Q (\Box_A^Q \, q') \leq q'$.  We multiply both sides of this by $\bigodot_Q\!\Gamma_Q$ and we get $f_A^Q (\Box_A^Q \, q') \bullet \bigodot_Q\!\Gamma \leq q' \bullet \bigodot_Q\!\Gamma$ and this is by the top line assumption less than $q$, so we have $f_A^Q (\Box_A^Q\, q') \bullet \bigodot_Q\!\Gamma \leq q' \bullet \bigodot_Q\!\Gamma \leq q$.
For the right rule our top line assumption is $f_A^Q (\bigodot_Q\!\Gamma) \leq q$ which is by adjunction equal to $\bigodot_Q\!\Gamma \leq \Box_A^Q \, q$. Note that this rule is also  sound on the other direction, that is the bottom line  implies the top line.

\smallskip
\noindent{\bf iii.~Soundness of epistemic update.} The rules for epistemic update  are
\[
\infer[\cdot  L]{m'\cdot q, \Gamma \vdash_M m}{m', q, \Gamma \vdash_M m}
\hspace{2cm}
\infer[\cdot  R]{\Gamma, \Gamma_Q, \Gamma_A \vdash_M m\cdot q}{\Gamma, \Gamma_A \vdash_M m \quad \Gamma_Q, \Gamma_A \vdash_Q q}
\] 
The soundness proofs for these rules use the definition of validity and satisfaction of $M$-sequents, which is based on the $\odot_M$ operation. So for the left rule we have to show the following
\[
\mbox{If} \quad m', q, \Gamma \models_M m \quad \mbox{then} \quad  m'\cdot q, \Gamma \models_M m
\]
which is by definition equivalent to the following 
\[
\mbox{If} \quad \bigodot_M (m', q, \Gamma) \leq m \quad \mbox{then} \quad  \bigodot_M (m'\cdot q, \Gamma) \leq m
\]
 This holds since  $\bigodot_M (m', q, \Gamma) =  \bigodot_M (m'\cdot q, \Gamma) = (m'\cdot q) \wedge \bigodot_M\!\Gamma$. 
 Proceeding similarly, for the right rule  we have to show the following
 \[
 \mbox{If} \quad \bigodot_M(\Gamma, \Gamma_A) \leq m \quad \mbox{and} \quad \bigodot_Q(\Gamma_Q, \Gamma_A) \leq q \quad \mbox{then} \quad \bigodot_M(\Gamma, \Gamma_Q, \Gamma_A) \leq m\cdot q
 \]
 In order to do so,   we first assume that  we have only one agent in our agent context, that is $\Gamma_A  = A$. By the first assumption of the top line we have  $f_A^M (\bigodot_M\!\Gamma) \leq m$ and by the second assumption we have $f_A^Q (\bigodot_Q\!\Gamma_Q) \leq q$. Since update is order preserving, we can update both sides of these two assumption by each other and get  $f_A^M (\bigodot_M\!\Gamma) \cdot  f_A^Q (\bigodot_Q\!\Gamma_Q) \leq m\cdot q$.  Now by update inequality we have $f_A^M (\bigodot_M\!\Gamma \cdot  \bigodot_Q\!\Gamma_Q) \leq f_A^M (\bigodot_M\!\Gamma) \cdot  f_A^Q (\bigodot_Q\!\Gamma_Q) \leq m\cdot q$, which is what we want for the bottom line and we are done. If we have more than one agent, that is $\Gamma_A = A_1, \dots,A_n$, then we follow the same line except that we have to apply the update inequality $n$ times, starting from the innermost agent $A_n$ to the outmost one $A_1$, that is 
\[
f_{A_n}^M (f_{A_{n-1}}^M (\dots f_{A_1}^M (\bigodot_M\!\Gamma\cdot \bigodot_Q\!\Gamma_Q))) \leq m\cdot q
\]

\smallskip
\noindent{\bf iv.~Soundness of  dynamic modality.}
The rules for dynamic modality are
\[
\infer[Dy L]{[q]m', \Gamma_Q \vdash_M m}{m' \vdash_M m \quad \Gamma_Q \vdash_Q q}
\hspace{2cm}
\infer[Dy R]{\Gamma \vdash_M [q]m}{\Gamma, q \vdash_M m}
\]
For the left rule we start from the second assumption $\bigodot_Q\Gamma_Q \leq  q$, since update is order preserving, we update $[q]m'$ on both sides  and we get
\[
[q]m' \cdot \bigodot_Q\!\Gamma_Q \leq [q]m' \cdot q
\]
Now by adjunction between update and dynamic modality we have  that $[q]m' \cdot q \leq m'$ and  by the first assumption of the top line we have $ m' \leq m$ and by transitivity we get
\[
[q]m' \cdot \bigodot_Q\!\Gamma_Q \leq  m
\]
which is exactly what we want for the bottom line.
We proceed similarly for the right rule, the $\odot_M$ definition of the  top line assumption says $\bigodot_M\!\Gamma\cdot q \leq m$, which is by adjunction equivalent to $\bigodot_M\!\Gamma \leq [q]m$ and  the $\odot_M$ definition of the  bottom line. Note that this rule holds in both directions, that is bottom line implies the top line.

\medskip
\begin{theorem}[Completeness] \em The rules presented in this section are complete with respect to the 
algebraic semantics in terms of distributive  epistemic systems.
\end{theorem}

\paragraph{Proof.} We show that if a sequent is valid in any distributive epistemic system  then it is provable using the rules of our $Q$ and $M$-systems. That is,  
\[
\mbox{if} \qquad \Gamma \models_Q q \qquad \mbox{then} \qquad \Gamma \vdash_Q q , \qquad \mbox{and} 
\qquad \mbox{if} \qquad  \Gamma \models_M m \qquad \mbox{then} \qquad \Gamma \vdash_M m\,.
\]
We show the contrapositive by building  two Lindenbaum-Tarski algebras: $M_0$ of equivalence classes of $M$-formulas over $\cong_M$   and $Q_0$ of equivalence classes of $Q$-formulas over $\cong_Q$ and define an order relation $\leq$ between them as $\vdash$ on their corresponding system. Similalry, we define all the algebraic operations of epistemic systems $\wedge, \vee, f_A, \Box_A, \cdot , [\ ], \bullet\,$ on the quantale and  module in terms of their sequent calculus counterparts,  and show that these operations  are well-defined over equivalence classes of formulas by using our sequent rules.  We then show that  these operations satisfy the   finite versions of  the equations of a distributive  
epistemic system. That is,  the same axioms but with  binary joins (and meets) instead of arbitrary ones.  Thus we have shown that $(M_0, Q_0, \{f_A\}_{A \in {\cal A}})$ constitutes a distributive  \emph{pre-epistemic system}, one with binary joins.  In order to extend our proof from this distributive pre-epistemic system to a distributive epistemic system (with arbitrary joins), we proceed by an  \emph{ideal construction}. We build the family of ideals over $M_0$ and $Q_0$, denoted by $M$ and $Q$,  and then show that  $(M, Q, \{f_A\}_{A \in {\cal A}})$ faithfully embeds  $(M_0, Q_0, \{f_A\}_{A \in {\cal A}})$ and thus it  is a   complete model of our sequent system. 

 The  full  proof is  presented  in~\cite{SadrThesis}, here  we proceed by providing the reader with some examples. In $Q_0$  the order is the logical consequence of $Q$-sequents   and the quantale operations are defined using the syntax of $Q$-formulas. Appearance maps and knowledge on $Q_0$ are defined using the $f_A^Q$ maps of the $Q$-system as follows
\[
f_A^Q ([q]) := [f_A^Q (q)]\qquad \mbox{and} \qquad  \Box_A^Q\,[q] = [\Box_A^Q\,q]
\]
We have to show that these are well-defined, that is
\[
\mbox{if} \quad [q_1] = [q'_1] \quad \mbox{then} \quad [f_A^Q(q_1)] = [f_A^Q(q'_1)]\,, \quad \mbox{and} \quad 
\mbox{if} \quad [q_1] = [q'_1] \quad \mbox{then} \quad [\Box_A^Q\,q_1] = [\Box_A^Q\,q'_1]
\]
or in logical consequence terms
\[
\mbox{if} \quad q_1 \vdash_Q\dashv q'_1 \quad \mbox{then} \quad f_A^Q(q_1)  \vdash_Q\dashv f_A^Q(q'_1)\,, \quad \mbox{and} \quad 
\mbox{if} \quad q_1  \vdash_Q\dashv q'_1 \quad \mbox{then} \quad \Box_A^Q\,q_1  \vdash_Q\dashv \Box_A^Q\,q'_1\,.
\]
 The proof trees for well-definedness of appearance are as follows
\[
\infer[f_A^Q L]{f_A^Q (q_1) \vdash_Q f_A^Q (q'_1)}{\infer[f_A^Q R]{q_1, A \vdash_Q f_A^Q (q'_1)}{\infer[Ass.]{q_1 \vdash_Q q'_1}{}}}
\hspace{2cm}
\infer[f_A^Q L]{f_A^Q (q'_1) \vdash_Q f_A^Q (q_1)}{\infer[f_A^Q R]{q'_1, A \vdash_Q f_A^Q (q_1)}{\infer[Ass.]{q'_1 \vdash_Q q_1}{}}}
\]
Similarly, the  proof trees for well-definedness of knowledge are
\[
\infer[\Box_A^Q R]{\Box_A^Q\,q_1 \vdash_Q \Box_A^Q\,q'_1}{\infer[\Box_A^Q L]{\Box_A^Q\, q_1, A \vdash q'_1}{\infer[Ass.]{q_1 \vdash_Q q'_1}{}}}
\hspace{2cm}
\infer[\Box_A^Q R]{\Box_A^Q\,q'_1 \vdash_Q \Box_A^Q\,q_1}{\infer[\Box_A^Q L]{\Box_A^Q\, q'_1, A \vdash q_1}{\infer[Ass.]{q'_1 \vdash_Q q_1}{}}}
\]
It remains to show that appearance and knowledge are adjoint, that is
\[
[f_A^Q (q)] \leq [q'] \quad \mbox{iff} \quad [q] \leq [\Box_A^Q\,q']
\]
The two proof trees for these follow
\[
\infer[f_A^Q L]{f_A^Q (q) \vdash_Q q'}{\infer[QCut]{q, A \vdash_Q q'}{\infer[Ass.]{ q \vdash_Q \Box_A^Q\,q' }{}& \infer[\Box_A^Q L]{\Box_A^Q\,q', A \vdash_Q q'}{\infer[Id]{q' \vdash_Q q'}{}}}}
\hspace{1.5cm}
\infer[\Box_A^Q R]{q \vdash_Q \Box_A^Q\,q'}{\infer[QCut]{q, A \vdash_Q q'}{\infer[f_A^Q R]{q, A \vdash_Q f_A^Q (q)}{\infer[Id]{q \vdash_Q q}{}} & \infer[Ass.]{f_A^Q (q) \vdash_Q q'}{}}}
\]
We now have to show that our operations satisfy the binary versions of axioms of epistemic systems. For example that the appearance maps on $Q_0$ preserve binary joins,  that is
\[
[f_A^Q (q_1 \vee q_2)] =  [f_A^Q (q_1) \vee f_A^Q (q_2)]
\]
The proof of the first direction of this equality is as follows
\[
\infer[f_A^Q  L]{f_A^Q (q_1 \vee q_2) \vdash_Q f_A^Q (q_1) \vee f_A^Q (q_2)}{\infer[\vee L]{q_1 \vee q_2, A \vdash_Q f_A^Q (q_1) \vee f_A^Q (q_2)}{
\infer[\vee R1]{q_1, A \vdash_Q f_A^Q (q_1) \vee f_A^Q (q_2)}{\infer[f_A^Q R]{q_1, A \vdash_Q f_A^Q (q_1)}{\infer[Id]{q_1 \vdash_Q q_1}{}}} \qquad
\infer[\vee R2]{q_2, A \vdash_Q f_A^Q (q_1) \vee f_A^Q (q_2)}{\infer[f_A^Q R]{q_2, A \vdash_Q f_A^Q (q_2)}{\infer[Id]{q_2 \vdash_Q q_2}{}}}
}}
\]
Similalry, the proof tree for the second direction is
\[
\infer[\vee L]{f_A^Q (q_1) \vee f_A^Q (q_2) \vdash_Q f_A^Q (q_1 \vee q_2)}{\infer[f_A^Q L]{f_A^Q (q_1) \vdash_Q f_A^Q (q_1 \vee q_2)}{\infer[f_A^Q R]{q_1, A \vdash_Q f_A^Q (q_1 \vee q_2)}{\infer[\vee R1]{q_1 \vdash_Q q_1 \vee q_2}{\infer[Id]{q_1 \vdash_Q q_1}{}}}}
\qquad
\infer[f_A^Q L]{f_A^Q (q_2) \vdash_Q f_A^Q (q_1 \vee q_2)}{\infer[f_A^Q R]{q_2, A \vdash_Q f_A^Q (q_1 \vee q_2)}{\infer[\vee R2]{q_2 \vdash_Q q_1 \vee q_2}{\infer[Id]{q_2 \vdash_Q q_2}{}}}}
}
\] 
The same constructions are done in the model built out of syntax of the $M$-system, that is  in $M_0$ where the order is $\vdash_M$. The meet, join, appearance and knowledge modality of $M_0$ are defined using their counterparts in the $M$-system, but for update and dynamic modality, we need to use both $M$ and $Q$-systems. The  update is defined on the pair $(M_0, Q_0)$ using the update operator of our $M$-systems as follows 
\[
[m] \cdot  [q] := [m \cdot  q]
\]
We have to show that it is well-defined, that is
\[
\mbox{If} \quad [m]= [m'] \quad \mbox{and} \quad [q]= [q'] \quad \mbox{then} \quad [m\cdot q] = [m'\cdot q']
\]
 The proof tree for one direction of this equality is as follows
\[
\infer[\cdot  L]{m \cdot  q \vdash_M m' \cdot  q'}{\infer[\cdot  R]{m,q \vdash_M m' \cdot  q'}{\infer[Ass.]{m \vdash_M m'}{} & \infer[Ass.]{q \vdash_Q q'}{}}}
\]
The proof tree for  the other direction is drawn similarly. It is easy to prove that update preserves binary joins of both $M_0$ and $Q_0$ and  the unit of $Q_0$, and that it is associative over multiplication of $Q_0$.  The dynamic modality of $M_0$ is defined in the same way by using the dynamic modality of the $M$-system and proved well-defined and adjoint to update. 

\smallskip
So far we have shown that $(M_0, Q_0, \{f_A\}_{A \in {\cal A}})$ is a distributive \emph{pre-epistemic} system with regard to which $M$ and $Q$-systems are complete. We extend this proof to distributive epistemic systems by  embedding this structure into an epistemic system  $(M,Q, \{f_A\}_{A \in {\cal A}})$ by taking $M = Idl(M_0)$ 
and $Q = Idl(Q_0)$ where  $Idl(M_0)$ is the family of \em ideals \em over $M_0$ and $Idl(Q_0)$ is the family of ideals over $Q_0$.  A subset of a lattice is called an 
\em ideal \em 
if it is non-empty, downward-closed, and closed under finite joins.
The order $\leq$ on ideals is given by inclusion, the arbitrary \em meet \em of ideals $\overline{\bigwedge_i} I_i$ is given 
by \em intersection \em of ideals $\bigcap_i I_i$, 
while the \em arbitrary join \em $\overline{\bigvee_i} I_i$ of a family of ideals is the ideal generated by their union, 
which is the downward-closure of the set of all finite joins of elements of these ideals. For example, the join in $Q$ is given by
\[
\overline{\bigvee_i} I_i :=\ \downarrow \!\left [ \left \{ \bigvee Y \Bigm | Y \ \mbox{finite} \ \subseteq \bigcup_{i} I_i \right \} \right ]\,.
\] 
The rest of operations, that is $f_A$ for both $M$ and $Q$, also $\cdot$ and  $\bullet\,$  are extended to ideals by applying them pointwise and then taking the  downward closure.  For instance, the appearance of ideals on $Q$ is defined as follows
\[
\overline{f_A^Q}(I) =\  \downarrow \!\!  [\{f_A^Q(q) \mid q \in I \}]
\]
We have to show that these operations are ideal preserving, that is for example, the join of ideals $\overline{\bigvee_i} I_i$ is an ideal.  Downward closure follows from the definition. For closure under joins assume that  $x,y \in \overline{\bigvee_i} I_i$, then $x \leq \bigvee Y_1$ and $y \leq \bigvee Y_2$, for $Y_1,Y_2$ finite subsets of the unions of $I_i$'s, that is $Y_1 \subseteq I_1$ and $Y_2 \subseteq I_2$. We have $x \vee y \leq (\bigvee Y_1) \vee  (\bigvee Y_2) = \bigvee (Y_1 \vee Y_2)$, since $Y_1\vee Y_2 \subseteq I_1 \cup I_2$, it is also a finite subset of union of $I_i$'s and thus $\bigvee(Y_1\vee Y_2)$ is an element of $\overline{\bigvee_i} I_i$. Since  $x \vee y$  lives in the  down set of $\bigvee(Y_1\vee Y_2)$,  we obtain $x \vee y \in \overline{\bigvee_i} I_i$. 
The proofs for other operations are done similarly, see~\cite{SadrThesis}.  The units of these operations are extended to ideals, the unit of multiplication is $\downarrow\!\!1$, the unit of appearance and join of $Q$ and $M$ is $\{\bot\}$ for the bottom of each accordingly, the unit of their meets is the ideal generated by the whole of $Q_0$ and $M_0$, that is $Q_0$ and $M_0$ themselves.  These ideals satisfy the axioms of epistemic systems, for example appearance of ideals of $Q_0$ preserves arbitrary joins of them. These are straightforward proofs and follow from the definition, for example for appearance of ideals of $Q_0$ we have to show
\[
\overline{f_A^Q} ( \overline{\bigvee_i} I_i) = \overline{\bigvee_i}\, \overline{f_A^Q} (I_i)
\]
We start from the left hand side
\begin{eqnarray*}
\overline{f_A^Q} ( \overline{\bigvee_i} I_i) &=& \downarrow\!\!\{f_A^Q(\overline{\bigvee_i} I_i) \mid I_i \ \mbox{is an ideal}\}\\
&=&  \downarrow\!\!\{f_A^Q(\bigvee Y)\mid Y \ \mbox{finite} \ \subseteq \bigcup_i I_i\}\\
&=&  \downarrow\!\!\{\bigvee f_A^Q(Y) \mid Y \ \mbox{finite} \ \subseteq \bigcup_i I_i\}
\end{eqnarray*}
which is equal to $\overline{\bigvee_i}\,\overline{f_A^Q} (I_i)$. The proofs for other axioms are done similarly and from them it follows that $(M,Q,\{f_A\}_{A \in {\cal A}})$ is an epistemic system. It remains to show that   $M_0$ and $Q_0$ are faithfully embedded into $M$ and $Q$. 
 The 
embedding $Q_0 \hookrightarrow Idl(Q_0)$ is defined as $q \mapsto \downarrow \!\! q$, and similarly for 
$M_0 \hookrightarrow Idl(M_0)$ as $m \mapsto \downarrow \!\! m$ . We show that this embedding  is a   homomorphism (thus it is faithful) in both $M$ and $Q$. We show this, for example in $Q$ and for $q_1, q_2 \in Q_0$,  by checking  the following
\begin{eqnarray*}
e(q_1)\,  \overline{\bigvee} \, e(q_2) &=& e(q_1\vee  q_2)\\
e(q_1)\,  \overline{\bigwedge} \, e(q_2) &=& e(q_1\wedge q_2)\\
e(q_1)\,  \overline{\bullet} \, e(q_2) &=& e(q_1 \bullet  q_2)\\
 \overline{f_A^Q} (e(q)) &=& e(f_A^Q (q)) 
\end{eqnarray*}
We   present the proof  for the  appearance maps of $Q$,  where we have to show $\overline{f_A^Q}(\downarrow \!\! q) \, = \,  \downarrow\!\!f_A^Q(q)$. By definition of appearance of ideals, this is equivalent to show  the following 
 \[
\downarrow\!\!\{f_A^Q (x) \mid \forall x \in \downarrow\!\!q\} \ = \ \downarrow\!\!f_A^Q (q)
\]
For the first direction, we take an element of the right hand side $
x \leq f_A^Q (q)$ 
and since $q \leq q$, we have $f_A^Q(q) \in \overline{f_A^Q} (\downarrow\!\!q)$  and we get $x \in  \overline{f_A^Q} (\downarrow\!\!q)$. 
For the other direction,  we take an element of the left hand side $x \in \overline{f_A^Q} (\downarrow\!\!q)$, which means  $x \leq f_A^Q(y)$ for some $y \leq q$. Since $f_A^Q$ is monotone, we apply it to both sides of $y \leq q$ and we get $f_A^Q(y) \leq f_A^Q(q)$, so we have that $
x \leq f_A^Q(q)$, that is an element of  the right hand side. 

\medskip
Since the distributive pre-epistemic system $(M_0,Q_0,\{f_A\}_{A \in {\cal A}})$ with binary operations was a complete model of our systems and the embedding is a homomorphism, we obtain that  the  distributive epistemic system $(M,Q,\{f_A\}_{A \in {\cal A}})$  inherits  completeness.   That is, we have the following  for the  $Q$-system   and a similar one for the $M$-system
\[
\text{If} \quad \Gamma_Q \nvdash_Q q \quad \text{then} \quad  e(\bigodot_Q \Gamma_Q) \nsubseteq_Q e(q)
\]
To see this, note that  from $\Gamma_Q \nvdash_Q q$  and completeness of the $Q_0$ system, it follows that $ \bigodot_Q \Gamma_Q \nleq_{Q_0} q$.  Since the embedding is a homomorphism, we have that   $  \downarrow\!\! \bigodot_Q\Gamma_Q \subseteq_Q  \downarrow\!\!q $ \ iff \ $\bigodot_Q \Gamma \leq_{Q_0} q$, from this  we get that  \  $ \bigodot_Q \Gamma_Q \nleq_{Q_0} q$ implies $e(\bigodot_Q \Gamma_Q) \nsubseteq_Q e(q)$. 
\endproof\newline

\bigskip
\noindent{\bf Example of  derivation}. \emph{Action-Knowledge Lemma and Prediction of Knowledge}. The lemma is as follows\footnote{It corresponds to a non-Boolean version  of the so-called
``Action-Knowledge Axiom'' of {\bf BMS}~\cite{BaltagMossSolecki}.}
\[
 \Box_A [f_A^Q(q)]m \vdash [q] \Box_A m\,.
\]
It uses an agent's  knowledge about the effect of his appearance of an action, that is $ \Box_A [f_A^Q(q)]m $ to derive his knowledge about the effect of  the action itself, that is $[q] \Box_A m$. It can also be seen as  a result about   permutation of epistemic $\Box_A$ and dynamic $[q]$ modalities up-to-appearance of actions $f_A^Q(q)$.

\noindent {\bf Proof.}
The proof tree is as follows

{\small
\[
\infer[Dy R]{\Box_A^M [f_A^Q(q)]m \vdash_M [q]\Box_A^M m}{
\infer[\Box_A^M R]{\Box_A^M [f_A^Q(q)]m, q \vdash_M \Box_A^M m}{
\infer[MCut]{\Box_A^M [f_A^Q(q)]m, q, A \vdash_M m}{
\infer[\cdot  R]{\Box_A^M [f_A^Q(q)]m, q, A \vdash_M [f_A^Q(q)]m \cdot  f_A^Q (q)}{
\infer[\Box_A^M L]{\Box_A^M[f_A^Q(q)]m, A\vdash_M [f_A^Q(q)]m}{
\infer[Id]{[f_A^Q(q)]m \vdash_M [f_A^Q(q)]m}{}}
&
\infer[f_A^Q R]{q, A \vdash_Q f_A^Q (q)}{\infer[Id]{q \vdash_Q q}{}}
}
&
\infer[\cdot  L]{[f_A^Q(q)]m \cdot  f_A^Q (q)\vdash_M m}{
\infer[Dy L]{[f_A^Q(q)]m, f_A^Q(q) \vdash_M m}{
\infer[Id]{m \vdash_M m}{}
&
\infer[Id]{f_A^Q (q) \vdash_Q f_A^Q (q)}{}
}}
}
}}
\]
}

\paragraph{Example of an application.} We  present the proof tree of the property that we proved for the MITM  cryptographic attack  in the algebra section. In order to encode the scenario in the sequent calculus, we have to add axioms for our appearances, facts, and kernel assumptions. For the appearance of propositions we have the following axiom schema for the $M$-systems (we refer to all these assumption axioms as $Ass.)$:
\[
\infer[Ass.]{m,A \vdash_M m'}{} \qquad \mbox{iff} \qquad f_A^M (m) = m' 
\]
Similarly, the following is  the axiom schema for the appearance of actions in the $Q$-system
\[
\infer[Ass.]{q,A \vdash_M q'}{} \qquad \mbox{iff} \qquad f_A^Q (q) = q' 
\]
For the kernel of actions, we have the following schema
\[
\infer[Ass.]{m,q \vdash_M \bot}{} \qquad \mbox{iff} \qquad m = Ker(q)
\]
and finally we encode the entailment between propositions  and facts $m \leq \phi$ via the following schema
\[
\infer[Ass.]{m \vdash_M \phi}{}\qquad \mbox{iff} \qquad m \leq \phi
\]
We encode the cryptographic attack scenario by instantiating these axioms. The axioms for the facts $P, \bar{P}$ and propositions $s,t$ will be the following 
\[
\infer[Ass.]{s \vdash_M P}{} \qquad \infer[Ass.]{t \vdash_M \bar{P}}{}
\]
We considered the kernel of  four actions $\{\alpha, \alpha', \beta, \beta'\}$ encoded as follows
\[
\infer[Ass.]{\bar{P},\alpha \vdash_M \bot}{} \quad \infer[Ass.]{\bar{P}, \alpha' \vdash_M \bot}{} \quad \infer[Ass.]{P, \beta \vdash_M \bot}{} \quad \infer[Ass.]{P, \beta' \vdash_M \bot}{}
\]
The encoding of the appearances of the propositions and actions to our three agents $\{A,B,C\}$  is straightforward, for example the ones used in the proof are encoded as follows on the $M$-system
\[
\infer[Ass.]{s,A \vdash_M s}{} \qquad \infer[Ass.]{s,B \vdash_M s\vee t}{}
\] 
and as follows for the actions in the $Q$-system
\[
\infer[Ass.]{\alpha,A \vdash_Q \alpha'}{}\qquad \infer[Ass.]{\alpha, B \vdash_Q \beta'}{}
\qquad \infer[Ass.]{\alpha', B \vdash_Q \alpha'}{}
\]
We prove that in the real state $s$ and after communicating the secret  $P$ or $\bar{P}$ via the action  $\alpha \vee \beta$, agent $A$ knows that $B$ knows that $P$ holds, that is  $s \cdot (\alpha \vee \beta) \vdash_M \Box_A^M \Box_B^M\,P$. One crucial part of the proof   is cut with an update formula and then the application of the left and right update rules to reduce the update to the assumptions axioms. The trick is to cut a  sequent that looks like $m,q,A \vdash_M m''$ with an update formula $m' \cdot  q'$ the proposition part of which is the appearance of the proposition on the left hand side, that is $f_A^M(m) = m'$ , and the action part of which is the appearance of the action on the left hand side, that is $f_A^Q(q) = q'$.  The other important part  of the proof is  cutting with $\bot$ and using the kernel assumption axioms. The steps of the proof are more or less the same as in the algebra. The proof tree  is as follows (in order to fit it in the page we had to draw two of its sub-trees $\mathbf \Pi1$ and $\mathbf \Pi2$ separately)

{\footnotesize
\[
\infer[\cdot  L]{s\cdot (\alpha \vee \beta) \vdash_M \Box_A^M \Box_B^M\,P}{
\infer[\vee ML]{s, \alpha \vee \beta \vdash_M  \Box_A^M \Box_B^M\,P}{
\infer[\Box_A^M R]{s, \alpha \vdash_M  \Box_A^M \Box_B^M\,P}{
\infer[MCut]{s, \alpha, A \vdash_M \Box_B^M\,P}{
\infer[\cdot R]{s, \alpha, A \vdash_M s\cdot \alpha'}{
\infer[Ass.]{s, A \vdash_M s}{}
&
\infer[Ass.]{\alpha, A \vdash_Q \alpha'}{}}
&
\infer[\cdot  L]{s\cdot \alpha' \vdash_M \Box_B^M\,P}{
\infer[\Box_B^M R]{s, \alpha' \vdash_M \Box_B^M\,P}{
\infer[MCut]{s, \alpha', B \vdash_M P}{
\infer[\cdot  R]{s, \alpha', B \vdash_M (s\vee t)\cdot \alpha'}{
\infer[Ass.]{s, B \vdash_M s\vee t}{}
&
\infer[Ass.]{\alpha', B \vdash_Q \alpha'}{}}
&
\infer[\cdot  L]{(s\vee t)\cdot \alpha' \vdash_M P}{\infer*{s \vee t, \alpha' \vdash_M P}{\mathbf \Pi1}
}}}}}}
&
\infer*{s, \beta \vdash_M  \Box_A^M \Box_B^M\,P}{\mathbf \Pi2}}}
\]
}
The sub-proof trees $\mathbf \Pi1$ and $\mathbf \Pi2$ are below

\bigskip\noindent
{\footnotesize
${\mathbf \Pi1}:$
\[
\infer[\vee L]{s \vee t, \alpha' \vdash_M P}{
\infer[fact]{s, \alpha' \vdash_M P}{
\infer[Ass.]{s \vdash_M P}{}}
&
\infer[MCut]{t, \alpha' \vdash_M P}{
\infer[MCut]{t, \alpha' \vdash_M \bot}{
\infer[Ass.]{t \vdash_M \bar{P}}{}
&
\infer[Ass.]{\bar{P}, \alpha' \vdash_M \bot}{}}
&
\infer[\bot L]{\bot \vdash_M P}{}}}
\]

\bigskip
\noindent
${\mathbf \Pi2:}$
\[
\infer[MCut]{s, \beta \vdash_M  \Box_A^M \Box_B^M\,P}{
\infer[MCut]{s, \beta \vdash_M \bot}{
\infer[Ass.]{s \vdash_M P}{}
&
\infer[Ass.]{P, \beta \vdash_M \bot}{}}
&
\infer[\bot L]{\bot \vdash_M \Box_A^M \Box_B^M\,P}{}}
\]
}

\section{Conclusion and further elaborations}     

We have developed an algebraic axiomatics in terms of a simple mathematical object: a sup-lattice $M$, 
which encodes  epistemic propositions and facts; a quantale $Q$ (acting on $M$) which 
encodes epistemic actions (and the updates induced by them); and a family of (lax-)endomorphisms of the  structure 
$(M,Q,\bigvee_M,\bigvee_Q,\cdot ,\bullet,1)$, encoding the agents' information states.  
From this structure many useful other modalities  arise, including dynamic 
modalities, epistemic modalities and residuals.  This algebraic axiomatics is a dynamic epistemic logic and generalizes the {\bf BMS} logic of~\cite{BaltagMossSolecki} to 
non-Boolean settings, while capturing the same concepts, 
and enriches it with a logical account of dynamic and epistemic resources in terms of actions and agents. We have presented  a sound and complete sequent calculus that enables us to deal with dynamic epistemic scenarios using  semi-automatic proof techniques. 
As examples of application, 
we have encoded and analyzed a classic epistemic puzzle (Muddy Children) and some of new variations of it with lying and 
cheating children, and proved the correctness of a simple security protocol, both algebraically and by a proof in 
sequent calculus. Some 
possible further elaborations on this line of thought follow.
\bit 
\item In this paper, following  dynamic epistemic logic, we dealt with the same update schema for all agents. This is a postulate of ``uniform rationality" and it means that the mechanism for  information update is the same for all agents.  It makes sense, if not being necessary, to consider personalized updates, where each agent updates his information in a different way than other agents do. We think that such personalized updates  could be better dealt with by moving to a categorical framework. More explicitly,  we are working in an enriched categorical setting where  a quantale $Q$ is a one-object quantaloid, i.e. a one-object sup-enriched category,  and agents' personalized updates ${\cal M}_A$ are sup-enriched functors. Appearance maps arise as lax sup-enriched natural transformations between the update functors. It would be interesting to compare our categorical approach  with  
coalgebraic epistemic features which are currently studied e.g.~\cite{Baltag3}. 
\item The Kripke semantics of a dynamic epistemic logic has been used as an alternative to BAN logic~\cite{BAN} to reason about security protocols e.g. in~\cite{Vink}. As shown in~\cite{SadrThesis} ch. 5, our algebraic setting provides an elegant frame work that facilitates these security applications. We would like to extend the domain of such applications  to be able to encode and prove the correctness of  open security protocols, for example by adding more types to our setting through a quantaloid enrichment~\cite{Stubbe}.
\item   {Approximation and probability.} We can conceive the modules in our setting as a more general 
type of partial orders than merely an algebraic logic. We can accommodate additional computational 
structure e.g.~a domain structure \cite{Compendium}, quantitative valuations of content \cite{Kopp, Martin}, or a 
combination of these which enables accommodating probabilities e.g.~the partial  
order on probability measures introduced in \cite{CoeMar} is defined in terms of a \em Bayesian update 
operation\em.  This development would also be of help in applications to security. 
\item Part of the motivation of this work was a marriage of epistemics and resource-sensitivity 
\cite{Mehrnoosh1}. Although we have introduced dynamic and epistemic resources in our setting, we would like to refine our logic and make it more resource-sensitive by relativizing our notion of ``consequence" to ``logical" actions available to agents. This will allow us to deal with classical   resource sensitive problems such as the problem of logical omniscience. The two examples below might provide useful 
insights, fragments and tools.  {\bf (i)}~~In the \em money games \em of \cite{Joyal} the resource, i.e., 
money 
$x\in\mathbb{R}^+$, is encoded using the quantale structure of $\mathbb{R}^+$ as a base for enrichment.  
The underlying lattices are free lattices which adds linearity to the propositions. They moreover admit a 
game-theoretic interpretation \cite{Joyal}.  {\bf (ii)}~~The \em logic of bunched implications \em of 
\cite{Pym} also provides a model to to handle resources which freely combines intuitionistic additive and 
multiplicative linear structure via contexts. The semantics in terms of Grothendieck sheaves 
of the additives again indicates a monoid-enriched structure in the sense of \cite{Stubbe}.
\item We would like to optimize our logic such that it has the cut-elimination property, this will involve  change of rules and might need a change of system. For example, and as suggested by our referee, an option would be to use the   deep inference deductive system in the \emph{calculus of structures}~\cite{Lutz}. We would also like to develop a boolean version of the sequent calculus presented here for concrete epistemic systems and prove its completeness with regard to Kripke semantics. Such a development will lead to a more refined version of our Theorem \ref{thm:represe} for a boolean dynamic epistemic logic.
\eit 
     
\section*{Acknowledgements} 
 
We thank Samson Abramsky, Andr\'e Joyal, Dusko Pavlovic,  Greg Restall, and 
Isar Stubbe for valuable discussions, and our referee Lutz  Stra{\ss}burger for his useful detailed  comments and corrections on the first version of this paper, and for his suggestion for the presentation of the sequent calculus. 
M.~S.~thanks Samson Abramsky and Oxford University Computer Laboratory for their hospitality, Mathieu Marion for his logistic support, and Roy Dyckhoff for his comments on cut-elimination. 
B.~C.~is supported by the EPSRC grant EP/C500032/1 High-Level Methods in Quantum Computation and Quantum Information.

\end{document}